\documentclass[11pt]{amsart}

\usepackage{latexsym, amsmath, amscd, graphicx}
\usepackage{epsfig, amssymb, amsfonts, mathrsfs, stmaryrd}
\usepackage[all]{xy}
\usepackage{hyperref}

\allowdisplaybreaks
\numberwithin{equation}{section}

\newcommand{\HH}{\mathbb{H}}
\newcommand{\Q}{\mathbb{Q}}
\newcommand{\R}{\mathbb{R}}

\def\Hom{{\rm Hom}}
\def\d{\delta}
\def\a{\alpha}
\def\b{\beta}

\def\D{\Delta}
\def\p{\partial}

\def\la{\longrightarrow}
\def\lm{\longmapsto}

\newtheorem{theorem}{Theorem}[section]
\newtheorem{lemma}[theorem]{Lemma}

\newtheorem{definition}[theorem]{Definition}
\newtheorem{def-thm}[theorem]{Definition-Theorem}

\newtheorem{corollary}[theorem]{Corollary}
\newtheorem{proposition}[theorem]{Proposition}

\theoremstyle{definition}

\newtheorem{remark}[theorem]{Remark}

\begin{document}

\title{An algebraic chain model of string topology}

\author{Xiaojun Chen}

\address{Department of Mathematics\\
University of Michigan\\
Ann Arbor, MI 48109}
\email{xch@umich.edu}

\date{December 26, 2008}

\maketitle

\begin{abstract}A chain complex model for the free loop space of a
connected, closed and oriented manifold is presented, and on its
homology, the Gerstenhaber and Batalin-Vilkovisky algebra structures
are defined and identified with the string topology structures. The
gravity algebra on the equivariant homology of the free loop space
is also modeled. The construction includes non simply-connected
case, and therefore gives an algebraic and chain level model of
Chas-Sullivan's String Topology.
\end{abstract}

\tableofcontents

\section{Introduction}\label{sectofintro}

This paper studies a chain complex
model of the free loop space of a smooth, compact and oriented
manifold. The purpose of our study is twofold: one is to give a
down-to-earth algebraic model of the structures: the
Gerstenhaber algebra, the Batalin-Vilkovisky algebra and the gravity
algebra, of string topology discovered by Chas-Sullivan in
\cite{CS99} and \cite{CS02} (see also Sullivan's survey \cite{Sullivan07}), and the other is to relate these algebraic
structures with some known ones, especially those from the
Hochschild complexes of the cochain algebra of the manifold.

The fundamental structure in our construction is the open Frobenius algebra of the chain complex of a
manifold. Given a smooth
compact manifold $M$, geometrically the chain complex $C_*(M)$ forms
a differential graded (DG) {\bf open Frobenius algebra}, namely, it is both a DG coalgebra
under the diagonal approximation and a DG algebra under the
transversal intersection, and the coproduct map is a map
of bimodules:
\begin{equation}\label{frobenius}\D(\a\cap\b)=\sum_{(\a)}\a'\otimes
\a''\cap\b=\sum_{(\b)}\a\cap\b'\otimes\b'',\end{equation}
where $\Delta$ is the diagonal map: $\Delta\a=\sum_{(\a)}\a'\otimes\a'',\Delta\b=\sum_{(\b)}\b'\otimes\b''$. This identity is called the Frobenius or module compatibility. However,
the Frobenius algebra is partially defined only when two chains are
transversal. For example, a chain of dimension less than that of the
manifold cannot intersect itself properly.

While most of the string topology operations are first defined on the chain level, it is conjectured that they may not be a homotopy invariant of the manifold (see \cite{Sullivan07}, Section 2.4 and also refer to \cite{CKS}),
a chain level intersection theory cannot be avoided. It turns out a weak version of the Frobenius algbra,
which we call the Frobenius-like algebra, on the Whitney forms of the manifold, is enough to model those
operations in string topology. Giving $M$ a smooth cubulation (by cubulation we mean a
decomposition of $M$ into cubes), recall that a Whitney polynomial
form on $M$ is a differential form on $M$ with rational polynomial coefficients
on each cube under the affine coordinates, with some obvious compatibilities. The set of Whitney polynomial forms, denoted by
$A(M)$, forms a DG algebra, whose homology is the rational
cohomology of $M$. Observe that by the compactness of $M$, the dual
space of $A(M)$, the set of currents, denoted by $C(M)$, forms a
complete DG coalgebra over $A(M)$. Moreover, $A(M)$ embeds into
$C(M)$ as in the smooth differential forms case, which is in fact a quasi-isomorphism by Poincar\'e duality.
If we view $A(M)$ as a subcomplex of $C(M)$,
and the wedge product on $A(M)$ as the intersection product (such a point of view is reasonable since on the homology level it does give the intersection product), then the induced
coproduct on $A(M)$: $$\D: A(M)\la C(M)\hat\otimes C(M)$$ satisfies the Frobenius compatibility
(\ref{frobenius}) formally (by ``formally" we mean the domain and range of $\Delta$ are in fact different). We would call the triple
$(A(M),C(M),\iota)$, where $\iota$ is the embedding map, a DG {\bf open Frobenius-like algebra}, whose homology
exactly gives the Frobenius algebra on $H_*(M)$ in the usual sense.

Note that $LM$ is a cosimplicial space (see Jones \cite{Jo87}), the associated cosimplicial chain complex, which is the complete cocyclic
cobar complex of $C(M)$, denoted by $\widehat{CC}_{*}(C(M))$, gives a chain complex model of $LM$ if $M$
is simply connected.
Here $\widehat{CC}_{*}(C(M))$ is the coalgebra analogue of the cyclic bar complex, and may also be viewed as the complete twisted tensor
product of $C(M)$ with its complete cobar construction $\hat\Omega(C(M))$ with a twisted differential (see below).
We call the homology of $\widehat{CC}_{*}(C(M))$
the {\bf Hochschild homology} of $C(M)$, and is denoted by $HH_{*}(C(M))$. Furthermore, the embedding of $A(M)$ into $C(M)$ together with Equation (\ref{frobenius}) guarantees that the linear map
$$\iota\otimes id: A(M)\hat\otimes\hat\Omega (C(M))\longrightarrow C(M)\hat\otimes \hat\Omega(C(M))=\widehat{CC}_{*}(C(M))$$
is a quasi-isomorphism of chain complexes (both with the twisted differential).

Both $A(M)$ and $\hat\Omega(C(M))$ are DG algebras; again by Equation (\ref{frobenius}) one can check that the twisted product
$A(M)\hat\otimes\hat\Omega(C(M))$ is a DG algebra under the twisted differential. The study of the commutativity
of the induced product leads to the following:

\begin{theorem}[Gerstenhaber algebra of the free loop space, cf. Theorem \ref{thmofgerst}]\label{thm1}
Given a DG open Frobenius-like algebra $(A, C,\iota)$, its Hochschild homology $HH_{*}(C)$ has the structure of a Gerstenhaber algebra. If the DG open
Frobenius-like algebra comes from a simply connected closed
manifold $M$, it gives the Gerstenhaber algebra on the homology of the free loop
space $LM$ with rational coefficients, which coincides with the one of Chas-Sullivan in string
topology.
\end{theorem}

Since the cocyclic cobar complex $\widehat{CC}_{*}(C(M))$ is the coalgebra analogue of the cyclic bar
complex, we may introduce the coalgebra analogue of the Connes cyclic operator
$$B:\widehat{CC}_{*}(C(M))\longrightarrow\widehat{ CC}_{*}(C(M)),$$
which models the $S^{1}$-rotation on $LM$ on the chain level. From
$A(M)\hat\otimes\hat\Omega(C(M))\simeq \widehat{CC}_{*}(C(M))$ we may pull back $B$ to the homology of $A\hat\otimes\hat\Omega(C(M))$
and obtain the following:

\begin{theorem}[Batalin-Vilkovisky algebra of the free loop space, cf. Theorem \ref{thmofbv}]
The functor in Theorem \ref{thm1} is in fact a functor to the category of
Batalin-Vilkovisky algebras, which gives the Batalin-Vilkovisky algebra of Chas-Sullivan in string topology in the case of simply connected manifolds.
\end{theorem}

The Batalin-Vilkovisky algebras are highly related to the topological field theories, see Getzler \cite{Ge94}. In fact, Getzler showed that a 2-dimensional (genus zero) topological conformal field theory (TCFT)
contains a natural Batalin-Vilkovisky structure. Later in \cite{Ge1} he continued to show that the equivariant TCFT (again in the genus zero case) has a structure of a
{\bf gravity algebra}, which may be viewed as a generalized Lie algebra.  By considering the cyclic
homology of the DG coalgebra $C(M)$, we obtain:

\begin{theorem}[Gravity algebra of the free loop space, cf. Theorem \ref{gravityalg}]
Given a DG open Frobenius-like algebra $(A,C,\iota)$, its cyclic homology has the structure of a gravity algebra, which models the gravity algebra of Chas-Sullivan on the equivariant
homology of the free loop space.
\end{theorem}

There is an extensive literature on string topology. Thomas Tradler, in his Ph.D. thesis \cite{Tr02}, first identified the loop homology (the homology of $LM$ with a degree shifting) with
the Hochschild cohomology of the the cochain complex of the manifold as Gerstenhaber algebras. His construction uses the singular chain complex of the manifold.
At the same time in the beautiful paper of Cohen-Jones \cite{CJ02} the authors gave a homotopy
theoretic realization of string topology via the Thom-Pontrjagin
construction, and also showed the isomorphism of the Hochschild cohomology with the loop homology. Voronov \cite{voronov1}
showed that the loop homology is an algebra over
the framed ``cactus" operad, while the latter is homotopic to the framed
little disk operad, which then implies the Batalin-Vilkovisky
algebra on the loop homology as well. For more details of their results, see
also Cohen-Hess-Voronov \cite{CHV}.

However, a shortcoming of the above approaches is that they are
mostly working on the homology level rather than on the chain level
(or at least assuming the manifold is formal). For example, it is
not easy to see from the homological construction of Cohen-Jones
that the loop product on the loop homology is commutative, which is
one of the key steps to the discovery of the Gerstenhaber algebra.
But this can be derived from our chain level construction (see Lemma
\ref{lemma4.3} of this paper). The chain level operations often
contain more of the structure of a manifold than those of homology.

McClure, in his paper \cite{Mc04}, constructed a chain level intersection theory by using the PL chains
(the intersection is partially defined), and from this he was able to show the Gerstenhaber and Batalin-Vilkovisky algebras on its cocyclic cobar complex.
In our construction of the Frobenius-like algebra, the set of currents is the chain complex of the manifold, and on its subcomplex, the Whitney forms, the intersection product is defined. Therefore
in some sense, the present paper may be viewed as a sequel to Tradler and McClure's work.

Our construction works over the field of rationals. For some other constructions using rational homotopy theory,
see F\'elix-Thomas-Vigu\'e \cite{FTV04},
Merkulov \cite{Me04} and F\'elix-Thomas
\cite{FT}.
Recently Menichi \cite{Menichi} gave a detailed analysis on the
Batalin-Vilkovisky algebra of string topology and of the Hochschild
cohomology of the cochain complex; in particular, he argued that if
the coefficient field is an arbitrary field, the situation is much
more delicate and complicated.

At last we remark that, according to the point of view of Dennis Sullivan, to correctly model the Frobenius structure of the
chain complex of a manifold,  one can: 1) either make the Frobenius compatibility Equation (\ref{frobenius}) strictly
hold but with the price that the intersection product is partially defined; 2) or diffuse the chains on the manifold such that the
intersection product is fully defined but with the price that Equation (\ref{frobenius}) only holds up to homotopy.
While the method applied in this paper takes the first point of view, the second point of view is also applicable: in the paper \cite{HL} Hamilton and Lazarev constructed a
cyclic Com$_{\infty}$ algebra which models the chain complex of a manifold, and then the result of
Costello (\cite{Costello}) shows that the Hochschild cohomology of the cyclic Com$_{\infty}$ algebras, which is isomorphic to the loop homology, naturally endows the structures of the Batalin-Vilkovisky algebra.
Recently McClure and Wilson (in private communication) have constructed a homotopy open Frobenius algebra on the chain
complex of a manifold, the application of their construction to model string topology, especially to obtain the higher genus string topology operations on the equivariant homology of the free loop space, will appear elsewhere.
A Thom-Pontrjagin method of these higher genus operations on the homology space, has been obtained by V. Godin (\cite{Godin}).

The rest of the paper is devoted to the proof of the above theorems. In
Section 2, we study the Frobenius algebra structure on the chain
complex of a smooth manifold. In Section 3 we construct a chain
complex model for the free loop space of a simply connected manifold. In Sections 4-6 we study the Gerstenhaber algebra
structure on its free loop space, from our point of
view. In Section 7 we
define and identify the Batalin-Vilkovisky algebra structure on the free loop
space. In Section 8, we recall a model for the equivariant chain
complex of the free loop space and identify the gravity algebra on its
homology. And in the last section we sketch the constructions of the
above structures on non-simply connected manifolds.

The author would like to thank Professor Dennis Sullivan, Professor
Yongbin Ruan, Professor James McClure, Professor John McCleary, Dr.
John Terilla, Dr. Thomas Tradler, Dr. Scott Wilson and Dezhen Xu for
many helpful conversations.

\section{The Chain Complex of a Manifold}

Let $M$ be a smooth, closed oriented manifold. The rational
cohomology of $M$, $H^*(M;\Q)$, has the following
structure: 1) it is a graded commutative algebra; 2) there is a
non-degenerate pairing on it by Poincar\'e duality. We usually call
a linear space which satisfies 1) and 2) a {\bf closed Frobenius
algebra}. Alternatively, a closed Frobenius algebra is a linear
space $V$ which is a graded commutative algebra and a graded
cocommutative coalgebra with both unit and counit, and moreover, the
coproduct is a map of bimodules:
\begin{equation}
\D(\a\cdot\b)=\sum_{(\a)}\a'\otimes\a''\cdot\b=\sum_{(\b)}\a\cdot\b'\otimes
\b''.\end{equation}By the isomorphism $H^*(M;\Q)\cong
H_{n+*}(M;\Q)$ (in this paper we grade the cohomology negatively, and the corresponding differential has degree $-1$),
the above statement says that $H_*(M;\Q)$ with the intersection
product and diagonal coproduct forms a closed Frobenius algebra.
However, such an algebraic structure does not hold on the chain
level, since the intersection of two chains is defined only when
they are transversal.

\begin{definition}[Whitney polynomial differential forms]\label{defofwhitney}
Let $M$ be a cubulated topological space. A cubical Whitney
polynomial differential form $x$ on $M$ is a collection of
differential forms, one on each cube, such that:
\begin{enumerate}\item[$(1)$] the coefficients of these forms on
each cube are $\Q$-polynomials with respect to the affine
coordinates of the cubes;
\item[$(2)$] they are compatible under restriction to faces, i.e. if
$\tau$ is face of $\sigma$, then
$x_\sigma|\tau=x_\tau$.\end{enumerate} The set of Whintney
polynomial forms on $M$ is denoted by $A(M)$.
\end{definition}

For a smooth manifold $M$ a smooth cubulation always
exists: by the famous theorem of Whitehead \cite{Wh}, any smooth
manifold admits a smooth triangulation, and therefore the dual
decomposition of such a triangulation naturally gives a smooth
cubulation of $M$. In the following we fix a smooth cubulation for
$M$.

Since $M$ is closed, the cubes on $M$ are finite in number, and
therefore if we denote by $A^p(M)$ the set of Whitney forms of
grading less than or equal to $p$ (here by grading we mean the sum
of the degree of the form and the degree of the polynomial
coefficient), then
\begin{equation}\label{filtration}A^0(M)\subset A^1(M)\subset\cdots,\quad \dim
A^p(M)<\infty\mbox{ for }p=0,1,\cdots\end{equation} and
$A(M)=\displaystyle\lim_{\longrightarrow}A^p(M)$. Moreover $A(M)$
has a unit and augmentation which are given by the constant
functions $A^0(M)\cong \Q$.

\begin{proposition}\label{Wh_fm}
Let $A(M)$ be the Whitney polynomial forms of $M$.
Then:\begin{enumerate}
\item[$(1)$]
$A(M)$, under wedge $\wedge$ and exterior differential $d$, forms a
commutative  DG algebra; \item[$(2)$] The Whitney forms may be
mapped to the cochains of the space as follows:
\begin{eqnarray*}
\rho: A(M)&\la&C^*(M;\Q)\\
x&\longmapsto& \Big\{I^n\mapsto\int_{I^n}x\Big\},\ \mbox{for any}\
I^n,
\end{eqnarray*}
which is a chain map.
\end{enumerate}
\end{proposition}

Proposition \ref{Wh_fm} (1) holds because $\wedge$ and $d$ are both natural
under restriction to faces, and (2) follows from Stokes' theorem.
Moreover, $A(M)$ computes the cohomology of $M$:

\begin{theorem}[de Rham's theorem for Whitney
forms]\label{theoremofderham} Let $M$ be a cubulated topological
space. Then $\rho$ is a chain equivalence of DG algebras, i.e.
$$\rho^*:H^*(A(M),d)\mathop{\la}^{\cong}_{\rm alg}H^*(M;\Q).$$
\end{theorem}

\begin{proof} See Cenkl-Porter \cite{CP}, Theorem 4.1.
\end{proof}

\begin{lemma}Let $M$ be a smooth manifold and $A(M)$ be the Whitney polynomial forms
of $M$. Let $C(M):=\Hom(A(M),\Q)$ be the space of currents; then $C(M)$ forms a differential
graded complete coalgebra with a counit and a coaugmentation.\end{lemma}

\begin{proof}
Note that
$C(M)=\Hom(A(M),\Q)=\Hom(\displaystyle\lim_{\la}A^p(M),\Q)=\displaystyle\lim_{\longleftarrow}\Hom(A^p(M),\Q)$,
and that the wedge product on $A(M)$ respects the filtration
(\ref{filtration}), $\wedge: A(M)\otimes A(M)\to A(M)$ induces a DG
map
\begin{eqnarray*}\D:C(M)=\Hom(A(M),\Q)&\longrightarrow& \Hom(A(M)\otimes
A(M),\Q)\\
&&=\Hom(\lim_{\la}A^p(M)\otimes\lim_{\la}A^p(M),\Q)\\
&&=\Hom(\lim_{\la}\bigoplus_{r=p+q}A^p(M)\otimes A^q(M),\Q)\\
&&=\lim_{\longleftarrow}\bigoplus_{r=p+q}\Hom(A^p(M),\Q)\otimes\Hom(A^q(M),\Q)\\
&&=C(M)\hat\otimes C(M),
\end{eqnarray*}
where the last equality holds by the definition of complete tensor
products: if $C=\displaystyle\lim_{\longleftarrow}C_p$ and
$D=\displaystyle\lim_{\longleftarrow}D_q$ are two inverse limit
systems, the complete tensor product of $C$ and $D$, denoted by
$C\hat\otimes D$, is given by
$$C\hat\otimes D:=\lim_{\longleftarrow}\bigoplus_{r=p+q}C_p\otimes
D_q.$$ The counit and coaugmentation come from the unit and
augmentation of $A(M)$.
\end{proof}

Since $A(M)$ computes the rational cohomology of $M$, by the
Universal Coefficient Theorem, $C(M)$ computes the rational homology
of $M$. We call $(C(M),\D, d)$ the {\bf complete DG
coalgebra} of $M$. As in the smooth case, the Whitney forms embed into the
currents, which is a quasi-isomorphism by Poincar\'e duality:

\begin{proposition}The embedding of $A(M)$ into $C(M)$, given
by
\begin{equation}
\label{dual}\iota:A(M)\la C(M): x\longmapsto\Big\{y\mapsto\int_M
x\wedge y\Big\},\end{equation} is a quasi-isomorphism of chain
complexes.
\end{proposition}

Note that $C(M)$ is a DG $A(M)$-bimodule, so if we denote $\D(\iota
x)$ by $\D x$, for any $x\in A(M)$, then:

\begin{proposition}
For any $x,y\in A(M)$,
\begin{equation}\label{froblike}\D(xy)=\sum_{(x)}x'\hat\otimes x''y=\sum_{(y)}xy'\hat\otimes
y''.\end{equation}\end{proposition}

The proof follows from a direct check. Note that Equation
(\ref{froblike}) is much like Equation (\ref{frobenius}), with
$\D:A(M)\to C(M)\hat\otimes C(M)$ instead of $\D:C(M)\to
C(M)\hat\otimes C(M)$. This shows that although we may not be able to define
the intersection and coproduct simultaneously on the chain complex
$C_*(M)$, the pair $(A(M),C(M))$ is good enough that we
may define the intersection on $A(M)$ and the coproduct on $C(M)$
while the Frobenius identity still holds.

Moreover, the coproduct of $A(M)$ factors through $A(M)\hat\otimes
C(M)$, namely, if we denote a basis of $A(M)$ by $\{y_i\}$, and let
$\{y_i^*\}$ be the dual basis, then $$\D
x=\sum_i\iota(xy_i)\hat\otimes y_i^*,$$ hence we may formally write
\begin{equation}\label{coprod_of_forms}
\D x=\sum_{i}xy_i\hat\otimes y_i^*,\quad\mbox{for all } x\in
A(M).\end{equation}

Let us summarize the above observations:
\begin{enumerate}\item[$(1)$] $A(M)$ is a DG commutative algebra;
\item[$(2)$] $C(M)$ is a (complete) DG
cocommutative coalgebra over $A(M)$;
\item[$(3)$] there is a quasi-isomorphic
embedding of $A(M)$-modules $\iota: A(M)\to C(M)$ which makes $A(M)$ a (complete) DG bi-comodule over $C(M)$.
\end{enumerate}

\begin{definition}
We call a triple $(A,C,\iota)$ which satisfies the above conditions $(1)$, $(2)$ and $(3)$
a DG {\bf open Frobenius-like algebra}.\end{definition}

The homology of $(A,C,\iota)$ is defined to be the homology of $A$
or $C$, which forms a Frobenius algebra naturally.
In the case of Whitney forms and currents on a manifold, this gives exactly
the closed Frobenius algebra structure on $H_*(M;\Q)$.

\section{Chain Complex Model of the Free Loop Space}

In this section we recall some facts on the cosimplicial chain complex model of the free loop space. The idea is due to K.-T. Chen \cite{Ch77} and Jones
\cite{Jo87} (see also Getzler-Jones-Petrack \cite{GJP91}).

Let $\D^n=\{(t_1,t_2,\cdots,t_n)\in\R^n:0\le t_1\le t_2\le\cdots\le
t_n\le 1\}$ be the standard $n$-simplex in $\R^n$. For each $n$, we
have the evaluation map
$$\Phi_n: LM\times \D^n\la\underbrace{M\times M\times\cdots\times M}_{n+1},$$
which is given by $\Phi_n(\gamma,
(t_1,t_2,\cdots,t_n))=(\gamma(0),\gamma(t_1),\cdots,\gamma(t_n))$.
By the chain equivalence of $C_*(M^{\times n+1})$ with
$C_*(M)^{\otimes n+1}$, consider the composition
$$C_*(LM)\la C_*(LM)\otimes [\D^n]\stackrel{\Phi_{n\#}}{\la} C_*(M)^{\otimes n+1},$$
where $[\D^n]$ is the fundamental chain of $\D^n$, $\Phi_{n\#}$ is the pushforward of $\Phi_n$, and we obtain a chain model for $LM$. Before doing that let us
introduce the {\bf cocyclic cobar complex} of a DG coalgebra, which
is the coalgebra analogue of the cyclic bar complex: Let $(C, d)$ be a DG
cocommutative coalgebra; the cocyclic cobar complex $CC_*(C)$ of $C$
is the direct product
$$\prod_{n=0}^\infty C\otimes (\Sigma C)^{\otimes n},$$
where $\Sigma$ is the desuspension functor of $C$ (the functor which
simply shifts the degree of $C$ down by 1), with differential
defined by
\begin{eqnarray}&& b(
a\otimes[a_1|\cdots|a_n])\label{twisteddiff}\\
&:=&da\otimes[a_1|\cdots|a_n]
-\sum_i(-1)^{|a|+|[a_1|\cdots|a_{i-1}]|}a\otimes [a_1|\cdots|da_i|\cdots|a_n]\label{difflpart}\\
&-&\sum_i\sum_{(a_i)}(-1)^{|a|+|[a_1|\cdots|a_{i-1}|a_i']|}a\otimes
[a_1|\cdots|a_i'|a_i''|\cdots|a_n]\label{diagonalpart}\\
&+&\sum_{(a)}(-1)^{|a'|}a'\otimes
\Big([a''|a_1|\cdots|a_n]-(-1)^{(|a''|-1)|[a_1|\cdots|a_n]|}[a_1|\cdots|a_n|a'']\Big).\label{diagonalpart_0}
\end{eqnarray}
Here we adopt the usual convention by writing elements in $C\otimes
(\Sigma C)^{\otimes n}$ in the form $a\otimes[a_1|\cdots|a_n]$. That
$b^2=0$ follows from the coassociativity of the coproduct of $C$.
Note that the dual complex of a DG coalgebra is a DG algebra; it is
easy to see that the dual complex of $CC_*(C)$ is the {\bf cyclic
bar complex} (also called the Hochschild complex) of the dual DG
algebra of $C$. We will call the homology of $CC_*(C)$ the {\bf
Hochschild homology} of the coalgebra $C$, denoted by $HH_*(C)$. For
an elegant treatment of the cyclic bar complex of a DG algebra, see,
for example, \cite{GJP91} and Loday \cite{loday}.

If moreover, $C$ is counital and coaugmented, we may consider the
{\bf reduced cocyclic cobar complex} of $C$, which is obtained from
$CC_*(C)$ by identifying elements $x\otimes [a_1|\cdots
|1|\cdots|a_n]$ with zero. To distinguish, we always write the latter
as $\prod_n C\otimes(\Sigma\bar C)^{\otimes n}$. In the
following when mentioning the cocyclic cobar complex we shall always
assume it is reduced, since in our construction of the DG coalgebra
of a manifold $M$, $C(M)$ is always counital and coaugmented.

We may extend the above definition to the case of complete DG
coalgebras, which is given by
$$\widehat{CC}_*(C):=\prod_{n=0}^{\infty}C\hat\otimes (\Sigma\bar C)^{\hat\otimes n},$$
with the differential $b$ extending to the complete tensor product.

\begin{theorem}\label{chain_equiv}Let $M$ be a simply connected manifold, and let $C(M)$ (written $C$ for short)
be the complete DG coalgebra of $M$. There is a chain equivalence
$$\phi:(C_*(LM),\partial)\la (\widehat{CC}_*(C),b).$$
\end{theorem}

\begin{proof}
The chain map is induced from \begin{eqnarray*}\phi: C_*(LM)&\la
&\displaystyle\prod_{n=0}^\infty C\hat\otimes (\Sigma\bar C)^{\hat\otimes n}\\
\a&\longmapsto&\displaystyle\sum \Phi_{n\#}(\a\otimes[\D^n]).
\end{eqnarray*}
Note that $\phi$ is a chain map: the differential of any element in
$\prod C\hat\otimes(\Sigma\bar C)^{\hat\otimes n}$ contains two
parts, one is those terms containing the differential of the
elements in $C$, the other is those terms that involve the coproduct
of the elements in $C$. If we write $b(\a)=b^{\rm I}(\a)+b^{\rm
II}(\a)$ referring to these two parts, namely, $b^{\rm
I}(\a)=(\ref{difflpart})$ and $b^{\rm
II}(\a)=(\ref{diagonalpart})+(\ref{diagonalpart_0})$, then
\begin{eqnarray*}
\phi(\p\a)&=&\sum \Phi_{n\#}\Big(\p\a\otimes[\D^n]\Big)\\
&=&\sum \Phi_{n\#}\Big(\p(\a\otimes[\D^n])-\a\otimes\p[\D^n]\Big)\\
&=&\p\circ\Big(\sum \Phi_{n\#}(\a\otimes[\D^n])
\Big)-\sum \Phi_{n\#}\Big(\sum_i(-1)^{i}\a\otimes\delta_i[\D^{n-1}]\Big)\\
&=&b^{\rm I}\circ\phi(\a)+b^{\rm II}\circ\phi(\a)=b\circ\phi(\a),
\end{eqnarray*}
where in the above $\delta_i$ is the identification of $\D^{n-1}$
with the $i$-th face of $\D^n$. More precisely, the last equality
holds due to the following: Define two groups of maps
$\{\delta_i:\Delta^{n-1}\to\Delta^n\}$ by
\begin{eqnarray*}\delta_0(t_1,\cdots,t_{n-1})&=&(0,t_1,\cdots,t_{n-1}),\\
\delta_i(t_1,\cdots,t_{n-1})&=&
(t_1,\cdots,t_i,t_i,\cdots,t_{n-1}),\quad 1\le i\le n-1\\
\delta_n(t_1,\cdots,t_{n-1})&=&(t_1,\cdots,t_{n-1},1),\end{eqnarray*}
and $\{\delta_i:M^{\times n}\to M^{\times n+1}\}$ by
\begin{eqnarray*}\delta_0(x_0,\cdots,x_{n-1})&=&(x_0,x_0,\cdots,x_{n-1}),\\
\delta_i(x_0,\cdots,x_{n-1})&=&
(x_0,\cdots,x_i,x_i,\cdots,x_{n-1}),\quad 1\le i\le n-1\\
\delta_n(x_0,\cdots,x_{n-1})&=&(x_0,\cdots,x_{n-1},x_0),\end{eqnarray*}
then the following diagram commutes
$$\xymatrix{
LM\times\Delta^{n-1}\ar[r]^-{\Phi_{n-1}}\ar[d]^{id\times\delta_i}&M^{\times n}\ar[d]^{\delta_i}\\
LM\times\Delta^n\ar[r]^-{\Phi_n}&M^{\times n+1} }$$ for all $0\le
i\le n$. Therefore,
$$\Phi_{n\#}\Big(\sum_i(-1)^{i}\a\otimes\delta_i[\D^{n-1}]\Big)=\sum_i(-1)^i\delta_{i\#}\circ\Phi_{n-1\#}(\a\otimes[\Delta^{n-1}]),$$
and if we shift the degree of the last $n-1$ components in
$C_*(M)^{\otimes n}$ down by one (which then multiplies $(-1)^i$ to
the image of $\delta_{i\#}\circ\Phi_{n-1\#}$, hence $(-1)^i$ cancel)
and sum over all $n\ge 0$, we obtain
$$\sum_{n=0}^\infty\Phi_{n\#}\Big(\sum_i(-1)^{i}\a\otimes\delta_i[\D^{n-1}]\Big)=b^{\rm
II}\circ\phi(\a).$$ The rest of the proof follows from a spectral
sequence argument; see, for example, Bousfield \cite{Bousfield},
\S4.1 or Rector \cite{Rector}, Corollary 5.2.
\end{proof}

We can also model the $S^1$-action on $LM$ in the above chain
complex model, given by the coalgebra version of {\bf Connes' cyclic
operator} (cf. Connes \cite{Connes} and Jones \cite{Jo87}):
$$\begin{array}{cccl}B:&\widehat{CC}_*(C)&\la&\widehat{CC}_*(C)\\
&a\otimes[a_1|\cdots|a_n]&\longmapsto&\displaystyle\sum_i(-1)^{|[a_i|\cdots|a_n]||[a_1|\cdots|a_{i-1}]|}
\varepsilon(a)a_i\otimes[a_{i+1}|\cdots|a_n|a_1|\cdots|a_{i-1}],
\end{array}$$ where $\varepsilon$ is the counit. One can easily check that $B^2=0$ and $bB+Bb=0$.

\begin{theorem}\label{chainequiv}
Let $M$ be a simply connected manifold. Let
$$J:C_*(LM)\la C_{*+1}(LM)$$
be the degree one map given by the composition
$$\begin{array}{ccccl}
LM&\la &LM \times S^1&\stackrel{A}{\la}&LM\\
\a&\longmapsto&\a \otimes[S^1]&\longmapsto&A_{\#}(\a\otimes[S^1]),
\end{array}$$ where $A$ is the rotation: $A(f,s)=f(s+\cdot)$, for any $f\in LM$ and $s\in S^1$, and $[S^1]$ is the fundamental cycle of $S^1$.
We have the following chain equivalence:
\begin{equation}(C_*(LM),d,J)\stackrel{\simeq}{\la}(\widehat{CC}_*(C),b,B).\label{Taction}\end{equation}
\end{theorem}

\begin{proof}Decompose $S^1\times \Delta^n$ into
$n+1$ standard $(n+1)$-simplices:
$S^1\times\Delta^n=\bigcup_{i=1}^{n+1}\Delta_i^{n+1}$, where
\begin{eqnarray}\D^{n+1}_i&:=&\{0\le
s\le\cdots\le s+t_{i-1}\le 1\le s+t_i\le\cdots\le s+t_n\le 2\}\nonumber\\
&=&\{0\le s+t_i-1\le\cdots\le s+t_n-1\le s\le\cdots\le s+t_{i-1}\le
1\},\label{changeofv2}
\end{eqnarray}
and let $r_i$ be the inclusion of $\Delta_i^{n+1}$ into
$S^1\times\Delta^n$, then the following diagram commutes
$$\xymatrix{
LM\times\Delta_i^{n+1}\ar@{^{(}->}[r]^{r_i}\ar[d]^{\Phi_{n+1}}&LM\times
S^1\times\Delta^n\ar[r]^-{A\times id}&LM\times\Delta^n\ar[d]^{\Phi_n}\\
M^{\times n+2}\ar[rr]^{\tau_i}&&M^{\times n+1}, }$$ where $\tau_i:
M^{\times n+2}\to M^{\times n+1}$ is given by
$$\tau_i(x_0,x_1,\cdots,x_{n+1})=(x_{n-i+1},\cdots,x_{n+1},x_1,\cdots,x_{n-i}),$$
and $\Phi_{n+1}$ is evaluated at
$(f,(0,s+t_i-1,\cdots,s+t_n-1,s,\cdots, s+t_{i-1}))$. Applying the
chain functor we obtain:
\begin{eqnarray*}\Phi_{n\#}(J\a\otimes[\Delta^n])&=&\Phi_{n\#}(A_{\#}(\a\otimes[S^1])\otimes[\Delta^n])\\
&=&\Phi_{n\#}\circ(A\times id)_{\#}(\a\otimes[S^1]\otimes[\Delta^n])\\
&=&\Phi_{n\#}\circ(A\times id)_{\#}\Big(\sum_{i=1}^{n+1} (id\otimes r_i)(\a\otimes[\Delta_i^{n+1}])\Big)\\
&=&\sum_{i=1}^{n+1}\tau_{i\#}\circ\Phi_{n+1\#}(\a\otimes[\Delta_i^{n+1}]).
\end{eqnarray*}
In particular, if
$\Phi_{n+1\#}(\a\otimes[\Delta^{n+1}])=a\otimes[a_1|\cdots|a_{n+1}]$,
then
$$\tau_{i\#}\circ\Phi_{n+1\#}(\a\otimes[\Delta^{n+1}_i])=\left\{\begin{array}{ll}
0,&\mbox{if}\quad |a|\ne 0;\\
\pm\displaystyle\varepsilon(a)
a_i\otimes[a_{i+1}|\cdots|a_{n+1}|a_1|\cdots|a_{i-1}],&\mbox{otherwise},
\end{array}\right.$$
where in the $|a|\ne 0$ case the value is zero because it is a
degenerate chain (the degrees of the two sides are not equal while
$\tau_{i\#}$ is a chain map), and therefore
\begin{eqnarray*}
\Phi_{n\#}(J\a\otimes[\D^n])&=&\sum_{i=1}^{n+1}\tau_{i\#}\circ\Phi_{n+1\#}(\a\otimes[\Delta_i^{n+1}])\\
&=&B\circ\Phi_{n+1\#}(\a\otimes[\D^{n+1}])
\end{eqnarray*}by definition. Summing over all $n\ge 0$, we obtain (\ref{Taction}) as claimed.
\end{proof}

In the above definition of $\widehat{CC}_*(C)$, if we write
$$\hat\Omega(C):=\prod_{n=0}^\infty (\Sigma\bar C)^{\hat\otimes n},$$
which is the {\bf complete cobar construction} of
$C$, then
$$\widehat{CC}_*(C)=C\hat\otimes\hat\Omega(C).$$
This has an interpretation of Brown's twisted tensor
product theory \cite{Br59}: For the fibration
$$\xymatrix{\Omega M\ar[r]&LM\ar[d]\\ &M,}$$
the theorem of Brown says that there is a chain equivalence between
the chain complex of the total space $LM$ and the ``twisted" tensor
product of the chain complexes of the base $M$ and the fiber $\Omega
M$. Since such a point of view plays a role in understanding the
Chas-Sullivan loop product on the loop homology, let us describe
this in more detail.

\begin{definition}[Twisting cochain]\label{defoftwistingcochain}
Let $(C,d)$ be a DG coalgebra over a field $k$ and $(A,\d)$ be a DG
algebra. A twisting cochain is a degree $-1$ linear map
$\Phi=\sum_q\Phi_q:C_q\to A_{q-1}$ such that
\begin{enumerate}\item[$(1)$] $\Phi_0(\varepsilon)=0$, where $\varepsilon$ is the counit;
\item[$(2)$] $\d\circ\Phi_q=-\Phi_{q-1} \circ
d-\sum_{k}(-1)^k\Phi_k\cup\Phi_{q-k}\circ\Delta$.
\end{enumerate}
\end{definition}

Let $(M,p)$ be a connected pointed topological space, and $S_*(M)$
be the 1-reduced singular chain complex of $M$ (here by ``1-reduced
singular chain complex" we mean the singular chains generated by
simplexes taking the vertices of the standard simplex into the
basepoint $p$). The Alexander-Whitney diagonal approximation gives a
DG coassociative coalgebra on $S_*(M)$. Now let $C_*(\Omega M)$ be
the chain complex of the based loop space of $M$ at the base point
$p$. Brown constructs a twisting cochain $\Phi:S_*(M)\to
C_{*-1}(\Omega M)$, which, roughly speaking, fills each simplex with
paths connecting its first and last vertices. Such a construction is
similar to the one of Adams' \cite{Ad56} (with a minor
modification). In fact, Adams proved that if $M$ is simply
connected, then the cobar construction of $S_*(M)$ is chain
equivalent to $C_*(\Omega M)$.

Now let $F\to E\stackrel{\pi}{\to} (M,p)$ be a Hurewicz fibration
with fiber $F=\pi^{-1}(p)$. Taking any loop $\gamma\in \Omega_{p}M$,
for any point $f\in F$ we may lift $\gamma$ to $E$ ending at $f$.
Denoting the initial point of the path to be $\gamma f$, we get a
continuous action of $\Omega_pM$ on $F$, which induces a DGA action
on the chain level:
$$\circ:\ C_*(\Omega M)\otimes C_*(F)\la C_*(F).$$
Define an operator $\partial_{\Phi}$ on $S_*(M)\otimes C_*(F)$ as
follows:
$$\partial_\Phi(a\otimes f):=\partial a\otimes f+(-1)^{|a|}a\otimes\partial f
+\sum(-1)^{|a'|}a'\otimes\Phi(a'')\circ f.$$ Then $\p_\Phi^2=0$. We
call $\p_\Phi$ the {\bf twisted differential} and $(S_*(M)\otimes
C_*(F),\partial_\Phi)$ the {\bf twisted tensor product}. The theorem
of Brown is that, for the fiber bundle $F\to E\to M$, there
is a chain equivalence
$$\phi:(S_*(M)\otimes
C_*(F),\partial_\Phi)\la(C_*(E),
\partial).$$

Now for the free loop space of a manifold $LM$, the fibration
$\Omega_p M\to LM\to (M,p)$ has a natural lifting function which is
given as follows: for any $\gamma: [0,1]\to M$, $\gamma(0)=q,\
\gamma(1)=p$, then
\begin{equation}\label{loopspaceaction}\begin{array}{rccl}
\gamma:&\Omega_pM&\la&\Omega_qM,\\
&x&\longmapsto&\gamma x\gamma^{-1}.
\end{array}\end{equation}
When $p=q$, the action is exactly the (left) adjoint action of
$\Omega_pM$ on itself. Passing to the chain level, it gives the
(left) adjoint action of the Hopf algebra $C_*(\Omega M)$ on itself.
Such an observation was also obtained by McCleary \cite{Mc88}.

By the result of Adams \cite{Ad56}, if $C$ is the DG coalgebra of
$M$, then the cobar construction $\Omega(C)$ gives a chain model
for $\Omega M$. It is not difficult to generalize the result to
the complete DG coalgebra case. Also the identity map $\tau:C\to
\hat\Omega(C): \a\mapsto [\a]$ is a twisting cochain, which exactly
models the one of Brown's. And therefore, the twisted tensor product
$$C\hat\otimes\hat\Omega(C)$$ with the twisted differential, which is the
differential $b$ given in Equation (\ref{twisteddiff}),
gives the chain complex model of $LM$.

Note that the Whitney polynomial forms $A(M)$ embed into $C(M)$;
thus we may form another chain complex
$$A\hat\otimes\hat\Omega(C)$$ with differential $b$ given by
\begin{eqnarray*}
&&b(x\otimes[a_1|\cdots|a_n])\\
&:=&dx\otimes[a_1|\cdots|a_n]-\sum_i(-1)^{|x|+|[a_1|\cdots|a_{i-1}]|}x\otimes
[a_1|\cdots|da_i|\cdots|a_n]\\
&-&\sum_i\sum_{(a_i)}(-1)^{|x|+|[a_1|\cdots|a_{i-1}|a_i']|}x\otimes[a_1|\cdots|a_i'|a_i''|\cdots|a_n] \\
&+&\sum_{i}(-1)^{|x|+|y_i|}xy_i\otimes
\Big([y_i^*|a_1|\cdots|a_n]-(-1)^{(|y_i|-1)|[a_1|\cdots|a_n]|}[a_1|\cdots|a_n|y_i^*]\Big).
\end{eqnarray*}
One can easily check that $b^2=0$. Comparing with Equation
(\ref{coprod_of_forms}), we see that it is also a twisted tensor
product, and by the comparison theorem of spectral sequences for
twisted tensor products (see e.g. McCleary \cite{McCleary} pp. 224), we have
that
\begin{equation}\label{restriction}\iota\hat\otimes id:
A\hat\otimes\hat\Omega(C)\la C\hat\otimes\hat\Omega(C)\end{equation}
is in fact a chain equivalence.

\section{The Chas-Sullivan Loop Product}

In this section we give a model of the Chas-Sullivan loop product
defined in string topology.

\begin{lemma}\label{modelofcs}
Let $(A,C,\iota)$ be a DG open Frobenius-like algebra. Define a product
$$
\bullet:A\hat\otimes\hat\Omega(C)\bigotimes
A\hat\otimes\hat\Omega(C)\longrightarrow A\hat\otimes\hat\Omega(C)$$
by
\begin{equation}\label{csloop}\big(x\otimes[a_1|\cdots|a_n]\big)\bullet
\big(y\otimes[b_1|\cdots|b_m]\big):=(-1)^{|y|)|[a_1|\cdots|a_n]|}x\wedge
y \otimes[a_1|\cdots|a_n|b_1|\cdots|b_m].\end{equation} Then
$(A\hat\otimes\hat\Omega(C),\bullet,b)$ forms a DG algebra.
\end{lemma}

\begin{proof}From the definition we see that $\bullet$ is
associative, so we only need to show $b$ is a derivation. Denoting
$x\otimes\a:=x\otimes[a_1|\cdots|a_n]$ and
$y\otimes\b:=y\otimes[b_1|\cdots|b_m]$ for short, up to sign, we
have
\begin{eqnarray}
&&b((x\otimes\a)\bullet(y\otimes\b))\label{b(abxy)}\\
&=&b(x y\otimes\a\cdot\b)\nonumber\\
&=&d(x y)\otimes\a\cdot\b+xy\otimes d(\a\cdot\b)\label{babxy}\\
&+&\sum(xy)'\otimes\tau (xy)''\circ(\a\cdot\b),\label{abtoxy}
\end{eqnarray}
where $\tau$ is the twisting cochain, which acts on $\hat\Omega(C)$
by the adjoint action, while
\begin{eqnarray}
&&b(x\otimes\a)\bullet (y\otimes\b)+(x\otimes \a)\bullet b(y\otimes
\b)\label{b(axby)}\\
&=& (dx) y\otimes \a\cdot\b+x y\otimes d(\a)\cdot\b\label{daxby}\\
&+&\sum x' y\otimes (\tau x''\circ \a)\cdot\b\label{atoxy}\\
&+&x(dy)\otimes\a\cdot\b+x y\otimes\a\cdot d(\b)\label{axdby}\\
&+&\sum x\cdot y'\otimes\a\cdot (\tau y''\circ\b).\label{btoax}
\end{eqnarray}
To show (\ref{b(abxy)})=(\ref{b(axby)}), noting that
(\ref{babxy})=(\ref{daxby})+(\ref{axdby}), we only need to show
(\ref{abtoxy})=(\ref{atoxy})+(\ref{btoax}), i.e.
$$\sum(xy)'\otimes\tau (xy)''\circ(\a\cdot\b)
=\sum x'y\otimes (\tau x''\circ\a)\cdot\b+\sum xy'\otimes\a\cdot
(\tau y''\circ\b).$$ By the Frobenius-like Equation (\ref{froblike})
it is equivalent for us to show
$$\tau z\circ(\a\cdot\b)=
(\tau z\circ\a)\cdot\b+\a\cdot( \tau z\circ\b),$$where $z=(xy)''$.
However, since all $\tau z$'s are primitive (for the definition and
properties of primitive elements of a complete Hopf algebra, see
e.g. Quillen \cite{Quillen}, Appendix A.2) and a direct calculation
shows that the primitive elements act as derivations, the result
follows.\end{proof}

Now let us briefly recall the {\bf loop product} defined in
\cite{CS99}. For the free loop space $LM$ of a manifold $M$, denote
by $C_*(LM)$ the chain complex of the total space. For an $x\in C_*(LM)$, suppose
it is not a linear combination of some other chains, then we may associate to $x$ a
chain $\tilde x\in C_*(M)$, which is the chain of marked points of $x$, and is called the ``shadow"
of $x$. Extend the map $x\mapsto
\tilde x$ linearly to all elements in $C_*(LM)$. Note that in general $x\mapsto\tilde x$ is not
a chain map. Now, for $x, y\in
C_*(LM)$ two chains in general position (transversal), the loop product of $x$ and $y$ is defined as follows:
first intersect $\tilde x$ and $\tilde y$ in $M$, then over the
intersection set, do the Pontrjagin product pointwisely. From this
we get a chain in $C_*(LM)$, denoted by $x\bullet y$, which is
usually called the {\bf loop product} of $x$ and $y$:
$$\begin{array}{cccl} \bullet:&C_*(LM)\otimes
C_*(LM)&\la&C_*(LM),\\
&x\otimes y&\lm& x\bullet y.
\end{array}
$$
Chas-Sullivan showed that $\partial$ is derivation with respect to
$\bullet$. A theorem of Wilson \cite{Wi04} says that although the
above product is defined on transversal chains, it already captures
all the homology information of $C_*(LM)$, and thus the loop product
is well-defined on the homology space $H_*(LM)$. Denote
$\HH_*(LM):=H_{*}(LM)[n]$ (which is called the {\bf loop homology}
of $M$); then $\HH_*(LM)$ is a graded algebra with the product
having degree 0.

\begin{theorem}[Model for the loop product]\label{modeloflp}
Let $M$ be a simply connected, smooth closed manifold. Then the
product $\bullet$ in Lemma \ref{modelofcs} gives a model of the loop
product in \cite{CS99}.
\end{theorem}

\begin{proof}Denote by
$$\phi:C_*(LM)[n]\la A\hat\otimes\hat\Omega(C)$$ the chain equivalence
(cf. Theorem~\ref{chain_equiv} and Equation (\ref{restriction})). In
the last section we have shown that $\phi$ is a chain map, so here
we only need to show $\phi$ is an algebra map. First let us consider
$\phi(x\otimes\a)$ and $\phi(y\otimes\b)$. They are two chains in
$LM$, whose geometric pictures are the chains swept by moving $\a$
(resp. $\b$) along $x$ (resp. $y$), for $x\otimes\a, y\otimes\b\in
A\hat\otimes\hat\Omega(C)$. Their shadows in $M$ are $x$ and $y$
respectively. Now $\phi(x\otimes \a)\bullet\phi(y\otimes\b)$ is a
chain in $LM$ described as follows: The shadow is $x\wedge y=xy$,
and for any point $q\in xy$, suppose there is a path $\gamma$
connecting $p$ and $q$, i.e.
$$\gamma:[0,1]\la xy\subset
M,\quad\gamma(0)=q,\quad\gamma(1)=p,$$ then by naturality of the
twisting cochain, the fiber over $q$ is the Pontrjagin product
\begin{equation}\label{fm}\gamma_{\#}(\a)\cdot\gamma_{\#}(\b),
\end{equation} where $\gamma_{\#}$ is
the chain map induced from
\begin{equation}\label{act}\begin{array}{cccl}
\gamma:&\Omega_pM&\la&\Omega_qM,\\
&\a&\lm&\gamma\cdot\a\cdot\gamma^{-1}. \end{array}\end{equation}

On the other hand, $\phi((-1)^{|\a||y|}xy\otimes\a\cdot\b)$ is a
chain in $LM$ described as follows: its shadow is also $xy$, and the
fiber over $q$ is
\begin{equation}\label{fp}
\gamma_{\#}(\a\cdot\b).
\end{equation}
In order to show
$$\phi(x\otimes\a)\bullet\phi(y\otimes \b)=\phi((-1)^{|\a||y|}xy\otimes\a\cdot\b),$$
we only need to show (\ref{fm})=(\ref{fp}):
\begin{equation}\label{eq}\gamma_{\#}(\a)\cdot\gamma_{\#}(\b)=\gamma_{\#}(\a\cdot \b).
\end{equation}
However, looking at the path action (\ref{act}), we have
$$\gamma(\a\cdot\b)=\gamma(\a)\cdot\gamma(\b),$$
for any $\a,\b\in\Omega_pM$, and on the chain level, it exactly
gives equality (\ref{eq}).
\end{proof}

\section{Commutativity of the Loop Product and the Gerstenhaber Algebra}

In \cite{CS99}, Chas and Sullivan showed that at the chain level,
the loop product is not commutative but homotopy commutative, and
hence the loop homology is commutative. Such a homotopy operator
gives a pre-Lie algebra on the loop homology, which was originally defined on the Hochschild cochain
complex of an associative algebra (see Gerstenhaber \cite{Ge63}).
The loop homology with the loop product and the commutator of the pre-Lie operator,
forms a Gerstenhaber algebra.

We first give a description of the pre-Lie operator $*$ defined in
\cite{CS99}: for two chains $\a,\b\in C_*(LM)$ in general position,
we have that $\tilde\a$ is transversal to loops in $\b$. Form a
chain $\a*\b$ given by the following loops: for any loop $\gamma$ in
$\b$, first go around $\gamma$ from the base point till the
intersection point with $\tilde\a$, then go around the loops in
$\a$, and finally go around the rest of $\gamma$. Such a star
$*$-operator can be modeled as follows:

\begin{lemma}\label{lemma4.3}
Let $(A,C,\iota)$ be as before. Define an operator
$$*:A\hat\otimes\hat\Omega(C)\bigotimes A\hat\otimes\hat\Omega(C)\la A\hat\otimes\hat\Omega(C)$$
as follows: for
$\a=x\otimes[a_1|\cdots|a_n],\b=y\otimes[b_1|\cdots|b_m]\in
A\hat\otimes\hat\Omega(C)$,
\begin{equation}\label{starop}\a*\b=\sum_{i=1}^n
(-1)^{|y|+(|\b|-1)|[a_{i+1}|\cdots|a_n]|}\varepsilon(a_i y)x\otimes
[a_1|\cdots|a_{i-1}|b_1|\cdots|b_m|a_{i+1}|\cdots|a_n],\end{equation}
where $\varepsilon$ is the counit of $C$. Then,
\begin{equation}\label{htpy}
b(\a*\b)=b\a*\b+(-1)^{|\a|+1}\a*b\b+(-1)^{|\a|}(\a\bullet\b-(-1)^{|\a||\b|}\b\bullet\a).
\end{equation}
In particular, $(H_*(A\hat\otimes\hat\Omega(C)),\bullet)$ is a graded
commutative algebra.
\end{lemma}

\begin{proof}The proof is more or less the same as Gerstenhaber \cite{Ge63}, Theorem
3, hence we omit it.
\end{proof}

\begin{definition}[Pre-Lie algebra]Let $V$ be a graded vector space over $k$.
A pre-Lie structure on $V$ is a degree one binary operator
$$*:V\otimes V\la V$$ such that
\begin{equation}\label{star}(\gamma*\a)*\b-
(-1)^{(|\a|+1)(|\b|+1)}(\gamma*\b)*\a=\gamma*(\a*\b-(-1)^{(|\a|+1)(|\b|+1)}\b*\a).
\end{equation}
We call $(V,*)$ a {\bf pre-Lie algebra}, or a {\bf pre-Lie system}.
\end{definition}

\begin{lemma}\label{gerprelie}
Let $(V,*)$ be a pre-Lie algebra. Define
$$\begin{array}{cccl}
\{,\}:&V\otimes V&\la& V\\
&a\otimes b&\longmapsto&a*b-(-1)^{(|a|+1)(|b|+1)}b*a,
\end{array}$$
then $(V,\{,\})$ is a degree one Lie algebra.
\end{lemma}

\begin{proof}
See Gerstenhaber \cite{Ge63}, Theorem 1.
\end{proof}

\begin{lemma}\label{prelieofst}
Let $(A,C,\iota)$ be as above. Then $(A\hat\otimes\hat\Omega(C),*)$
is a pre-Lie algebra.
\end{lemma}

\begin{proof}We also omit the proof; one may refer to Gerstenhaber \cite{Ge63}, Sections 5-7.\end{proof}

\begin{corollary}\label{degreeonelie}
Let $(A,C,\iota)$ be as above. Then
$$(A\hat\otimes\hat\Omega(C),\{,\},b)$$ is a degree one DG Lie
algebra. In particular, $(H_*(A\hat\otimes\hat\Omega(C)),\{,\})$ is a
degree one graded Lie algebra.
\end{corollary}

\begin{proof}
The degree one Lie algebra follows from the above lemma and the
theorem of Gerstenhaber (Lemma \ref{gerprelie}). Lemma
\ref{lemma4.3} shows that $b$ respects $\{,\}$.
\end{proof}

\begin{definition}[Gerstenhaber algebra]\label{defofgerst}
Let $V$ be a graded vector space over a field $k$. A Gerstenhaber
algebra on $V$ is a triple $(V,\cdot,\{,\})$ such that
\begin{enumerate}\item[$(1)$] $(V,\cdot)$ is a graded commutative
algebra; \item[$(2)$] $(V,\{,\})$ is a graded degree one Lie
algebra; \item[$(3)$] the bracket is a derivation for both
variables.
\end{enumerate}
\end{definition}

We are now ready to model the theorem of \cite{CS99}, where the Lie
bracket $\{,\}$ is called the {\bf loop bracket:}

\begin{theorem}[Gerstenhaber algebra of the free loop sapce]
\label{thmofgerst}Let $M$ be a simply connected, smooth closed
manifold and $LM$ its free loop space. Then
$$(H_*(A\hat\otimes\hat\Omega(C)),\bullet,\{,\})$$
is a Gerstenhaber algebra, which models the Gerstenhaber algebra on
$\HH_*(LM)$ obtained in \cite{CS99}.
\end{theorem}

\begin{proof}We have shown that $H_*(A\hat\otimes\hat\Omega(C)$ is a
graded commutative algebra (Lemma \ref{lemma4.3}) and a degree one
graded Lie algebra (Corollary \ref{degreeonelie}). Next we show that
the bracket is a derivation with respect to the loop product for
both variables. By symmetry we only need to show, for
$\a,\b,\gamma\in H_*(A\hat\otimes\hat\Omega(C))$,
$$\{\a\bullet\b,\gamma\}=\a\bullet\{\b,\gamma\}+(-1)^{|\b|(|\gamma|+1)}
\{\a,\gamma\}\bullet\b.$$ This immediately follows from the
following Lemma \ref{lemma4.5}.

As we shall see later (Section \ref{sectofBV}), both $\HH_*(LM)$ and
$H_*(A\hat\otimes\Omega(C))$ has a Batalin-Vilkovisky algebra
structure, and $\{,\}$ is completely determined by the
Batalin-Vilkovisky differential operator. The identification of the
Batalin-Vilkovisky algebras (Theorem \ref{thmofbv}) then gives the
identification of the Gerstenhaber algebras of $\HH_*(LM)$ and
$H_*(A\hat\otimes\hat\Omega(C))$. However, the identification of
$\{,\}$ cannot be obtained directly from the above arguments, even
though we followed Gerstenhaber \cite{Ge63} and Chas-Sullivan
\cite{CS99} step by step, since it comes from the commutator of the
homotopy operator $*$, which is {\it a priori} not a chain map.
\end{proof}

\begin{lemma}\label{lemma4.5}Let $A$ be as above.
Then for
$\a=x\otimes[a_1|\cdots|a_n],\b=y\otimes[b_1|\cdots|b_m],\gamma=z\otimes
[c_1|\cdots|c_l]\in A\hat\otimes\hat\Omega(C)$,
\begin{enumerate}\item[$(1)$]
$(\a\bullet\b)*\gamma=\a\bullet(\b*\gamma)+(-1)^{|\b|(|\gamma|+1)}(\a*\gamma)
\bullet\b;$ \item[$(2)$] setting
\begin{eqnarray*}&&h(\gamma\otimes\a\otimes\b)\\
&=&\sum_{i<j}(-1)^{\epsilon}\varepsilon(c_ix)\varepsilon(c_jy)
z\otimes[c_1|\cdots|c_{i-1}|a_1|\cdots|a_n|c_{i+1}|\cdots
|c_{j-1}|b_1|\cdots|b_m|c_{j+1}|\cdots|c_l],
\end{eqnarray*}
where
$\epsilon=(|\a|-1)|[c_{i+1}|\cdots|c_l]|+(|\b|-1)|[c_{j+1}|\cdots|c_l]|$,
we have
$$(b\circ h-h\circ b)(\gamma\otimes\a\otimes\b)
=\gamma*(\a\bullet\b)-(\gamma*\a)\bullet\b-(-1)^{|\a|(|\gamma|+1)}\a\bullet(\gamma*\b).$$
\end{enumerate}
\end{lemma}

The above lemma is similar to \cite{CS99}, Lemma 4.6, with a
minor modification, whose proof is deferred to the Appendix.

\section{Isomorphism of Two Gerstenhaber Algebras}

The notion of a Gerstenhaber algebras was first
introduced by Gerstenhaber in his study of the deformation theory of
associative algebras. In \cite{Ge63} Gerstenhaber showed that the
Hochschild cohomology of an associative algebra is endowed with the
structure of a Gerstenhaber algebra. As a direct application, the
Hochschild cohomology of the cochain algebra $C^*(M)$ of a manifold
is a Gerstenhaber algebra. As we have seen, the (co)homology of the
free loop space is closely related to the appropriate Hochschild
homology of the (co)chain algebra; one wonders if the Gerstenhaber
algebra deduced from Gerstenhaber's result is identical to the one
discovered in string topology.

Such a problem has been discussed and proved by Cohen-Jones \cite{CJ02},
Tradler \cite{Tr02}, Merkulov \cite{Me04}, F\'elix-Thomas-Vigu\'e
\cite{FTV04} and McClure \cite{Mc04}. Here we also offer an affirmative
answer via our chain model of the free loop space.
Recall the results of Gerstenhaber in \cite{Ge63}:

\begin{definition}[Product and bracket of the Hochschild cochain complex]
Let $A$ be a (DG) algebra over a field $k$ and let
$$HC^*(A;A)=\Hom\big(\bigoplus_{n=0}^\infty A^{\otimes n},A\big)$$ be its Hochschild cochain complex. Define the
product $\cup$, the pre-Lie operator $*$, and the bracket $\{,\}$ on
$HC^*(A;A)$ as follows: for $f\in\Hom(A^{\otimes n};A)$,
$g\in\Hom(A^{\otimes m};A)$, up to sign,
\begin{enumerate}\item[$(1)$]$f\cup g\in\Hom(A^{\otimes m+n},A)$: for any $a_1,\cdots,a_{m+n}\in A$,\begin{equation}
\label{gerprod}(f\cup g)(a_1,\cdots,a_{m+n}):=f(a_1,\cdots,a_n)\cdot
g(a_{n+1},\cdots,a_{m+n});\end{equation}
\item[$(2)$]$f*g\in\Hom(A^{\otimes m+n-1},A)$: for any $a_1,\cdots,a_{n+m-1}\in A$,
\begin{equation}\label{gerstar}
(f*g)(a_1,\cdots,a_{n+m-1}):=\sum_{i=1}^{n}\pm
f(a_1,\cdots,a_{i-1},g(a_i,\cdots,a_{i+m-1}),\cdots,a_{n+m-1});
\end{equation}
\item[$(3)$] $\{f,g\}\in\Hom(A^{\otimes m+n-1},A)$ is the commutator of $*$:
\begin{equation}\label{gerbracket}
\{f,g\}:=f*g-(-1)^{(|f|+1)(|g|+1)}g*f.
\end{equation}
\end{enumerate}
\end{definition}

Gerstenhaber's main statement in \cite{Ge63} is the following
theorem:

\begin{theorem}[Gerstenhaber]
Let $A$ be a DG associative algebra over a field $k$ and let the
operators $\cup$, $*$ and $\{,\}$ be given in the above definition;
then Lemmas (\ref{lemma4.3}) and (\ref{lemma4.5}) hold. Therefore
the Hochschild cohomology $(HH^*(A;A),\cup,\{,\})$ is a Gerstenhaber
algebra.
\end{theorem}

The following theorem is inspired by the results of the authors
mentioned at the beginning of this section:

\begin{theorem}[Isomorphism of two Gerstenhaber algebras]\label{thmofiso} Let $M$ be
a simply connected manifold and $A$ be the Whitney forms on $M$.
Then
$$\HH_*(LM)\stackrel{\cong}{\la}HH^*(A;A)$$
are isomorphic as Gerstenhaber algebras.
\end{theorem}

\begin{proof}
In fact, let $C$ be the set of currents on $M$; then the Hochschild
cochain complex is chain equivalent to $A\hat\otimes\hat\Omega(C)$:
$$HC^*(A;A)\simeq A\hat\otimes\hat\Omega(C).$$
For homogeneous $f,g\in HC^*(A;A)$, we may write them as
$f=x\otimes[a_1|\cdots|a_n],g=y\otimes[b_1|\cdots|b_m]\in
A\hat\otimes\hat\Omega(C)$; the operators $\cdot$, $*$ and $\{,\}$
defined above by (\ref{gerprod}), (\ref{gerstar}) and
(\ref{gerbracket}) can be rewritten as
\begin{equation*}f\cup g=x\cdot y\otimes[a_1|\cdots|a_n|b_1|\cdots|b_m]
\end{equation*}
and
\begin{equation*}
f*g=\sum_{i=1}^n \langle a_i,y\rangle
x\otimes[a_1|\cdots|a_{i-1}|b_1|\cdots|b_m|a_{i+1}|\cdots|a_n],
\end{equation*}
and $$\{f,g\}:=f*g-(-1)^{(|f|+1)(|g|+1)}g*f.$$ Comparing them with
the loop product (\ref{csloop}) and pre-Lie operator (\ref{starop}),
we see that $\HH_*(LM)$ and $HH^*(A;A)$ are isomorphic as
Gerstenhaber algebras.
\end{proof}

\begin{remark}In the paper of Voronov and Gerstenhaber \cite{VG}, the authors show that the Hochschild cochain complex has a very ample structure, which
they called the {\bf homotopy Gerstenhaber algebra}, where a family
of new operators besides the pre-Lie operator are introduced: while
the pre-Lie operator gives the homotopy of the commutativity, these
new operators give all the higher homotopies. From the above proof
we may see that the chain complex of the free loop space of a
manifold also bears the structure of a homotopy Gerstenhaber
algebra, which is highly related to the cactus operad and the little
disk operad (see Voronov \cite{voronov} and \cite{voronov1}).
\end{remark}

\section{$S^1$-action and the Batalin-Vilkovisky
Algebra}\label{sectofBV}

Let $J$ be the $S^1$-action on the loop homology. In \cite{CS99},
Chas and Sullivan prove that $(\HH_*(LM),\bullet,J)$ forms a
Batalin-Vilkovisky algebra. Namely, $J$ on homology is not a
derivation with respect to $\bullet$, but the deviation from being a
derivation is a derivation. One deduces that,
$$\{a,b\}:=(-1)^{|\a|}J(\a\bullet\b)-(-1)^{|\a|}J(\a)\bullet b-\a\bullet J(\b),\quad\mbox{for all}\quad
\a,\b\in\HH_*(LM),$$
defines a degree one graded Lie algebra on $\HH_*(LM)$, which is in
fact the loop bracket on homology.

\begin{definition}[Batalin-Vilkovisky
algebra]\label{batalinvilkovisky} Let $V$ be a graded vector space
over a field $k$. A Batalin-Vilkovisky algebra on $V$ is a triple
$(V,\bullet,\D)$ such that:
\begin{enumerate}\item[$(1)$] $(V,\bullet)$ is a graded commutative
algebra; \item[$(2)$] $\D:V\to V$ is degree one operator with
$\D^2=0$; \item[$(3)$]The deviation from being a derivation of $\D$
with respect to $\bullet$ is a derivation for both variables,
namely,
$$(-1)^{|\a|}\D(\a\bullet\b)-(-1)^{|\a|}\D(\a)\bullet b-\a\bullet
\D(\b)$$ is a derivation for both $\a,\b\in V$.
\end{enumerate}
\end{definition}

A Batalin-Vilkovisky algebra is a special kind of Gerstenhaber
algebra:

\begin{proposition}
Let $(V,\bullet,\D)$ be a Batalin-Vilkovisky algebra. Define $[\
,]:V\otimes V\la V$ by
$$[\a,\b]:=(-1)^{|\a|}\D(\a\bullet\b)-(-1)^{|\a|}\D(\a)\bullet
b-\a\bullet \D(\b), \mbox{ for } \a,\b\in V,$$ then $(V,\bullet,[\
,])$ forms a Gerstenhaber algebra.
\end{proposition}

\begin{proof}See Getzler \cite{Ge94}, Proposition 1.2. \end{proof}

\begin{lemma}\label{bvgersten}
Let $M$ be a simply connected manifold and $LM$ be its free loop
space. Then
\begin{equation}\label{bv}
\{\a,\b\}=(-1)^{|\a|}J(\a\bullet\b)-(-1)^{|\a|}J(\a)\bullet
b-\a\bullet J(\b),\mbox{ for }\a,\b\in\HH_*(LM),
\end{equation}
where $\{,\}$ and $\bullet$ are the loop bracket and the loop
product respectively, and $J$ is the induced $S^1$-action on
$\HH_*(LM)$.

More precisely, let $(A, C,\iota)$ be the DG open Frobenius-like
algebra of $M$ and $A\hat\otimes\hat\Omega(C)$ be the twisted tensor
product, and let $B$ be the dual of Connes' cyclic operator on
$C\hat\otimes\hat\Omega(C)$. Then there is a linear map
$$h:A\hat\otimes\hat\Omega(C)\bigotimes A\hat\otimes\hat\Omega(C)\la
C\hat\otimes\hat\Omega(C)$$ such that for any $\a,\b\in
A\hat\otimes\hat\Omega(C)$,
\begin{equation}\label{bvhtpy}
(b\circ h-h\circ b)(\a\otimes
\b)=\{\a,\b\}-(-1)^{|\a|}B(\a\bullet\b)-(-1)^{(|\b|+1)(|\a|+1)}\b\bullet
B(\a) +\a\bullet B(\b).
\end{equation}
\end{lemma}

The above lemma is similar to Lemma 5.2 in \cite{CS99}, whose proof is also given in the Appendix. By this
lemma we obtain:

\begin{theorem}[Batalin-Vilkovisky algebra of the free loop space]
\label{thmofbv} Let $M$ be a simply connected, smooth closed
manifold and let $A$ be the Whitney forms and $C$ be the currents on
$M$. Then
$$(H_*(A\hat\otimes\hat\Omega(C)),\bullet,B)$$
is a Batalin-Vilkovisky algebra, which models the Batalin-Vilkovisky
algebra on $\HH_*(LM)$ obtained in \cite{CS99}.
\end{theorem}

\begin{proof}We have shown (Theorem \ref{thmofgerst}) that
$$(H_*(A\hat\otimes\hat\Omega(C)),\bullet,\{,\})$$
is a Gerstenhaber algebra, and therefore the loop bracket $\{,\}$ is
a derivation for both variables with respect to $\bullet$. Lemma
\ref{bvgersten} says that the deviation of $B$ from being a
derivation is exactly the loop bracket. Thus, according to
Definition \ref{batalinvilkovisky},
$$(H_*(A\hat\otimes\hat\Omega(C),\bullet,B)$$ is a
Batalin-Vilkovisky algebra.

In Theorem~\ref{chainequiv} we identified $H_*(LM)$ with
$H_*(\widehat{CC}_*(C),b)$ and hence
$H_*(A\hat\otimes\hat\Omega(C),b)$ (up to a degree shifting) as
vector spaces, together with the identification of the
$S^1$-rotation with Connes' cyclic operator $B$. In
Theorem~\ref{modeloflp} we identified the loop product with the
product on $H_*(A\hat\otimes\hat\Omega(C)$. From the definition, a
Batalin-Vilkovisky algebra is completely determined by these two
operations, and therefore the Batalin-Vilkovisky algebra obtained
above models the one of string topology.
\end{proof}

\section{Equivariant Homology and the Gravity Algebra}

In Chas-Sullivan \cite{CS99} the authors also showed that the
equivariant homology of the free loop space, $H_*^{S^1}(LM)$, forms
a Lie algebra. Later in \cite{CS02} they continued to show that the
equivariant homology is endowed with more structures, for example, the gravity algebra. Recall that the equivariant homology
$H_*^{S^1}(LM)$ of $LM$ is defined to be the homology of
$ES^1\times_{S^1}LM$, where $ES^1$ is the universal bundle over the
classifying space $BS^1$. There is a fibration
$$\xymatrix{S^1\ar[r]&ES^1\times LM\ar[d]\\
&ES^1\times_{S^1}LM,}$$ the associated Gysin sequence is given by
$$\cdots\la H_*(ES^1\times LM)\la H_*^{S^1}(LM)\la
H_{*-2}^{S^1}(LM)\la H_{*-1}(ES^1\times LM)\la\cdots.$$ By
identifying $H_*(ES^1\times LM)$ with $H_*(LM)$ we obtain
$$\cdots\la H_*(LM)
\stackrel{E}{\la} H_*^{S^1}(LM)\la
H_{*-2}^{S^1}(LM)\stackrel{M}{\la} H_{*-1}(LM)\la\cdots,$$ where $E$
and $M$ can be interpreted as ``forgetting" and ``adding" the marked
points on the loops of the corresponding spaces.

Topologically $M\circ E$ is exactly the $S^1$-operation $J$ on
homology $H_*(LM)$, and $E\circ M=0$. Now for $a_1,a_2\in
H_*^{S^1}(LM)$, define
$$\{a_1,a_2\}:=(-1)^{|a_1|+2-n}E(M(a_1)\bullet  M(a_2)),$$
which is usually called the {\bf string bracket}, then $\{,\}$ thus
defined in fact gives on $H_*^{S^1}(LM)$ a degree $2-n$ graded Lie
algebra structure. Moreover, $H_*^{S^1}(LM)$ satisfies
the generalized Jacobi identity, and hence forms a gravity algebra
in the sense of Getzler \cite{Ge1}:

\begin{definition}[Gravity algebra]Let $V$ be a chain complex over
a field $k$. A gravity algebra on $V$ is a sequence of graded
skew-symmetric operators:
$$c_n:V^{\otimes n}\la V,\quad\mbox{for}\,\,n\ge 2,$$
of degree $2-n$, satisfying the following relations: if $k>2$ and
$l\ge 0$, and $a_1,\cdots,a_k, b_1,\cdots, b_l\in V$,
\begin{eqnarray}
&&\sum_{1\le i< j\le k}(-1)^\epsilon
\{\{a_i,a_j\},a_1,\cdots,\widehat{a_i},\cdots,\widehat{a_j},
\cdots,a_k,b_1,\cdots,b_l\}\label{identity_of_grav}\\
&=&\left\{\begin{array}{lc}
\{\{a_1,\cdots,a_k\},b_1,\cdots,b_l\},&l>0,\\
0,&l=0,
\end{array}\right.\nonumber\end{eqnarray}
where we write $c_n(a_1,\cdots, a_n)$ as $\{a_1,\cdots, a_n\}$, and
$\epsilon=|a_i|(|a_i|+\cdots+|a_{i-1}|)+|a_j|(|a_1|+\cdots+|\widehat{
a_i}|+\cdots+|a_{j-1}|)$.\end{definition}

A gravity algebra contains a Lie algebra: let $k=3$ and $l=0$; then
Equation (\ref{identity_of_grav}) is exactly the Jacobi identity.
For more details of the gravity algebra on the equivariant homology
$H_*^{S^1}(LM)$, see \cite{CS99}, \cite{CS02}, \cite{Su03} or
Theorem \ref{gravityalg} below. The above construction is rather
topological, but we can see this algebraically from the cyclic
homology of A. Connes.

\begin{definition}[Cyclic homology of a coalgebra]Let $C$ be a DG
coalgebra. The cyclic homology of $C$, denoted by $CH_*(C)$, is the
homology of the chain complex $CC_*(C)[u,u^{-1}]/\langle u^{-1}\rangle$, where $u$
is a parameter of degree $2$, with differential induced from the one
defined on $CC_*(C)[u,u^{-1}]$:
$$\begin{array}{cccl}b+u^{-1}B:&CC_*(C)[u,u^{-1}]&\la &CC_*(C)[u,u^{-1}]\\
&x\otimes u^n&\longmapsto& bx\otimes u^n+ Bx\otimes u^{n-1}.
\end{array}$$\end{definition}

According to Jones \cite{Jo87}, there are several cyclic homology
theories. Here we adopt the most common used one in literature. The above definition can be generalized to the complete DG
coalgebra case.

\begin{theorem}[Connes' exact sequence and the Gysin sequence]
$(1)$ Let $C$ be a DG cocommutative coalgebra. Then there is a long
exact sequence, called Connes' exact sequence:
\begin{equation}\label{Connes_exact_seq}
\xymatrix{\cdots\ar[r]&HH_*(C)\ar[r]^E&CH_*(C)\ar[r]&CH_{*-2}(C)\ar[r]^M&HH_{*-1}(C)\ar[r]&\cdots.}
\end{equation}

$(2)$ If moreover, $C$ is the DG coalgebra of a simply connected
manifold $M$, then the following diagram is commutative:
$$\xymatrix{
\cdots\ar[r]&H_*(LM)\ar[d]^{\cong}\ar[r]&H_*^{S^1}(LM)\ar[d]^{\cong}
\ar[r]&H_{*-2}^{S^1}(LM)\ar[d]^{\cong}\ar[r]&H_{*-1}(LM)\ar[d]^{\cong}\ar[r]&\cdots\\
\cdots\ar[r]&HH_*(C)\ar[r]^{E}&CH_*(C)\ar[r]&CH_{*-2}(C)\ar[r]^{M}&HH_{*-1}(C)\ar[r]&\cdots\\
}$$
\end{theorem}

\begin{proof}The proof of the two statements is the coalgebra
analogue of Loday \cite{loday}, Theorem 7.2.3, p. 235. In fact,
observe that we have a short exact sequence:
$$0\la CC_*(C)\la
CC_*(C)[u,u^{-1}]/\langle u^{-1}\rangle\stackrel{\cdot
u^{-1}}{\la}CC_*(C)[u,u^{-1}]/\langle u^{-1}\rangle\la 0.$$ The associated long
exact sequence on homology is exactly Connes' long exact sequence.
The isomorphism between
$$H_*^{S^1}(LM)\stackrel{\cong}{\la}CH_*(C)$$
comes from the fact that $C_*^{S^1}(LM)$ is chain equivalent to (see
Jones \cite{Jo87})
$$(C_*(LM)[u,u^{-1}]/\langle u^{-1}\rangle,b+u^{-1}J).$$
Applying Theorem \ref{chainequiv} gives the desired isomorphism.
\end{proof}

\begin{lemma}In the long exact sequence (\ref{Connes_exact_seq}) of
the above theorem,
$$M\circ E=B:HH_*(C)\la HH_{*+1}(C).$$
\end{lemma}
\begin{proof}
The statement follows from chasing the diagram of the short exact
sequence
$$0\la \widehat{CC}_*(C)\la \widehat{CC}[u,u^{-1}]/u^{-1}\stackrel{\cdot
u^{-1}}{\la}\widehat{CC}_*(C)[u,u^{-1}]/u^{-1}\la 0.$$ The check is
left to the reader.\end{proof}

\begin{theorem}[Gravity algebra on the free loop
space]\label{gravityalg} Let $M$ be a simply connected manifold and
let $C$ be the DG coalgebra of $M$. Let
$\mathbb{CH}_*(C):=CH_*(C)[n-2]$, and define
\begin{eqnarray*}
c_n:\mathbb{CH}_*(C)^{\otimes n}&\la&\mathbb{CH}_*(C)\\
\a_1\otimes\cdots\otimes\a_n&\longmapsto&(-1)^{\epsilon}E(M(\a_1)\bullet\cdots\bullet
M(\a_n)),
\end{eqnarray*}
where $\bullet$ is the loop product, and
$\epsilon=(n-1)|\a_1|+(n-2)|\a_2|+\cdots+|\a_{n-1}|$. Then
$(\mathbb{CH}_*(C),\{c_n\})$ is a gravity algebra.
\end{theorem}

\begin{proof}We have shown that $(HH_*(C),\bullet,B)$ is a
Batalin-Vilkovisky algebra. $B$ being a second order operator with
respect to $\bullet$ implies that \begin{eqnarray}B(x_1\bullet
x_2\bullet\cdots\bullet x_n)&=&\sum_{i<j}\pm B(x_i\bullet x_j)\bullet
x_1\bullet\cdots\bullet\widehat{x_i}\bullet\cdots\bullet\widehat{x_j}
\bullet\cdots\bullet x_n\label{deriv}\\
&\pm& (n-2)\sum_i x_1\bullet\cdots\bullet B x_i\bullet\cdots\bullet
x_n,\nonumber\end{eqnarray} for $x_1,\cdots, x_n\in HH_*(C)$.

Now let $x_i:=M(a_i)$, and apply $E$ on both sides of the above
equality to obtain:
\begin{eqnarray*}&&E\circ B(M(a_1)\bullet
M(a_2)\bullet\cdots\bullet M(a_n))\\
&=&\sum_{i<j}\pm E\circ\Big( B(M(a_i)\bullet M(a_j))\bullet
M(a_1)\bullet\cdots\bullet\widehat{M(a_i)}\bullet\cdots\bullet\widehat{M(a_j)}
\bullet\cdots\bullet M(a_n)\Big)\\
&\pm& (n-2)\sum_i E\circ\Big(M(a_1)\bullet\cdots\bullet B\circ
M(a_i)\bullet\cdots\bullet M(a_n)\Big).\end{eqnarray*} Note that
$E\circ B=E\circ M\circ E=0$ and $B\circ M=M\circ E\circ M=0$ (above
lemma), so we exactly have $$\sum_{1\le i< j\le k}\pm
\{\{a_i,a_j\},a_1,\cdots,\widehat{a_i},\cdots,\widehat{a_j},
\cdots,a_k\}=0.$$ Similarly by multiplying $y_1\bullet\cdots\bullet
y_l$ to (\ref{deriv}), letting $y_j:=M(b_j)$ and applying $E$ on
both sides, we obtain
$$\sum_{1\le i< j\le k}\pm
\{\{a_i,a_j\},a_1,\cdots,\widehat{a_i},\cdots,\widehat{a_j},
\cdots,a_k,b_1,\cdots,b_l\}=\{\{a_1,\cdots,a_k\},b_1,\cdots,b_l\},$$
for $l>0$. This proves the theorem.
\end{proof}

\section{The Non-simply Connected Manifolds}

In the previous sections, we have only discussed the case when the
manifold $M$ is simply connected. In this section we sketch the
construction of string topology on a general non-simply connected
manifold. The idea is to lift the loops on $M$ to its universal
covering $\tilde M$, where the loops now becomes paths, which can be
characterized explicitly. This idea is due to Mike Mandell, which
was communicated to the author by James McClure. Since $\tilde M$ is
simply connected, our algebraic methods may now be applied.

We begin with the following observation about the free loop space
$LM$.

\begin{lemma}[Equivalent characterization of $LM$]\label{equivchar}
Let $M$ be a smooth manifold. Denote by $G$ the fundamental group
$\pi_1(M)$ and by $\tilde M$ the universal covering of $M$. For any
$g\in G$, let
$$L_g\tilde M:=\Big\{f:I=[0,1]\to\tilde M \Big|f(1)=g\circ
f(0)\Big\}.$$ Then $\coprod_{g\in G}L_g\tilde M$ admits a $G$-action
induced from that on $\tilde M$: for $f\in L_g\tilde M$, and $h\in
G$,
$$\begin{array}{cccc}h\circ f:&[0,1]&\la&\tilde M\\
&x&\longmapsto &h\circ f(x).
\end{array}$$
Since $(h\circ f)(1)=h\circ f(1)=h\circ (g\circ
f(0))=hgh^{-1}\circ((h\circ f)(0))$, $h\circ f\in L_{hgh^{-1}}\tilde
M$. There is a homeomorphism
$$\coprod_{g\in G}L_g\tilde M/G\cong LM,$$
and the following diagram commutes:
$$\xymatrixcolsep{4pc}\xymatrix{
\displaystyle\coprod_{g\in G} L_g\tilde M\ar[r]^{/G}\ar[d]^{\pi_0}&LM\ar[d]^{\pi_0}\\
\tilde M\ar[r]^{/G}&M, }$$ where $\pi_0$ is the projection of the
paths to their starting points.
\end{lemma}

The proof of the lemma is a direct check. In the following we shall use this proposition to construct a
chain complex model for $LM$ from $\coprod_gL_g\tilde M$. However,
since $\tilde M$ may not be closed, the dual space of the Whitney
polynomial differential forms may not compute the homology of
$\tilde M$ correctly, and therefore may not be a correct chain model
for $\tilde M$. However, if we denote by $A_{\sigma}^p$ the set of
Whitney forms of degree less than or equal to $p$ on a cube $\sigma$
in $\tilde M$, then by the definition of the Whitney forms,
$$A(\tilde M)={\lim_{\longleftarrow}}_\sigma{\lim_{\longrightarrow
}}_pA_{\sigma}^p.$$ Let $C_{\sigma}^p:=\Hom(A_\sigma^p,\Q)$, and
$$C(\tilde M):={\lim_{\longrightarrow}}_\sigma{\lim_{\longleftarrow}}_{p}C_\sigma^p.$$
$C(\tilde M)$ may be viewed as the set of currents with compact
support. Similar to the cochain case, there is a chain map
$$\rho:C_*(\tilde M)\la C(\tilde M)$$ from the singular chain complex $C_*(\tilde M)$ to
$C(\tilde M)$ which is given by integration, inducing an
isomorphism on the homology. Denote by $A_{\rm c}(\tilde M)$ the set
of Whitney forms with compact support. Then there is an embedding
$$\begin{array}{cccl}
\iota:&A_{\rm c}(\tilde M)&\la&C(\tilde M)\\
&\a&\longmapsto&\displaystyle\left\{\b\mapsto\int_{\tilde
M}\a\wedge\b\right\}.\end{array}$$ By the fact that $H_{\rm
c}^*(\tilde M)\cong H_*(\tilde M;\Q)$ one deduces that $$( A_{\rm
c}(\tilde M), C(\tilde M),\iota)$$ is a DG open Frobenius algebra.
Moreover, the action of $G$ on $\tilde M$ induces a $G$-action on
$A_{\rm c}(\tilde M)$ and $C(M)$, and the inclusion $\iota: A_{\rm
c}(\tilde M)\la C(\tilde M)$ is in fact $G$-equivariant.

Recall the definition of $L_g\tilde M$:
$$L_g\tilde M=\Big\{f:[0,1]\to\tilde M \Big|f(1)=g\circ
f(0)\Big\}.$$ Note that $L_g\tilde M$ is connected: given
$f_1,f_2\in L_g\tilde M$, let $\gamma$ be a path in $\tilde M$
connecting $f_1(0)$ and $f_2(0)$. Then $g\circ\gamma$ is a path
connecting $f_1(1)$ and $f_2(1)$, and $f_1\circ (g\circ\gamma)\circ
f_2^{-1}\circ \gamma^{-1}$ is a closed path in $\tilde M$. Since
$\tilde M$ is simply connected, $f_1\circ (g\circ\gamma)\circ
f_2^{-1}\circ \gamma^{-1}$ can be filled in with paths, which gives
a path in $L_g\tilde M$ connecting $f_1$ and $f_2$.

Now consider the evaluation maps (compare with \S3)
\begin{equation}\label{evaluationmaps}
{\Psi_n}:L_g\tilde M\times\Delta^n {\la}\tilde
M\times\cdots\times\tilde M
\end{equation}
given by
$$\Psi_n(f,(t_1,\cdots,t_n)):=(f(0),f(t_1),\cdots,f(t_n)).$$
On the chain level this leads to a chain complex which is rather
similar to the cocyclic cobar complex:

\begin{definition} Let $(C,\D,d)$ be a
coassociative DG coalgebra over field $k$. Suppose $G$ is a discrete
group and $C$ admits a $k[G]$-action, which commutes with $\D$. Let
$\Omega(C)$ be the cobar construction of $C$. Fix $g\in G$. Define
an operator
$$b_g:C\otimes\Omega(C)\la C\otimes\Omega(C)$$
by
\begin{eqnarray*}
&&b_g(x\otimes[a_1|\cdots|a_n])\\
&:=&dx\otimes[a_1|\cdots|a_n]-\sum_i(-1)^{|x|+|[a_1|\cdots|a_{i-1}]|}x\otimes[a_1|\cdots|da_i|\cdots|a_n]\\
&-&\sum_i\sum_{(a_i)}(-1)^{|x|+|[a_1|\cdots|a_{i-1}|a_i']|}x\otimes[a_1|\cdots|a_i'|a_i''|\cdots|a_n]\\
&+&\sum_{(x)}(-1)^{|x'|}\Big(x'\otimes[x''|a_1|\cdots|a_n]-
(-1)^{(|x''|-1)(|[a_1|\cdots|a_n]|)}x'\otimes[a_1|\cdots|a_n|g_*x'']\Big),
\end{eqnarray*}then $b_g^2=0$.
\end{definition}

Consider the direct sum of $(C\otimes\Omega(C), b_g)$ indexed by
$G$, and denote it by
$$\big(C\otimes\Omega(C)
\otimes k[G],\tilde b=\sum_{g\in G}b_g\big).
$$ And define a $k[G]$-action on it by $$\begin{array}{ccl} \left(\displaystyle k[G],
C\otimes\Omega(C) \otimes k[G]\right)&\la&C\otimes\Omega(C)
\otimes k[G]\\
(h,x\otimes[a_1|\cdots|a_n]\otimes g)&\longmapsto&
h_*x\otimes[h_*a_1|\cdots|h_*a_n]\otimes hgh^{-1}.
\end{array}$$
Moreover, define an operator $\tilde B$ on $C \otimes
\Omega(C)\otimes k[G]$ as follows:
$$\displaystyle
\begin{array}{rccl}
\tilde B:& C \otimes\Omega(C)\otimes k[G]&\longrightarrow&
C\otimes\Omega (C)\otimes k[G]\\
&x\otimes[a_1|\cdots|a_n]\otimes
g&\longmapsto&\displaystyle\sum_i(-1)^{\epsilon}
\varepsilon(x)a_i\otimes[a_{i+1}|\cdots|a_n|g_*a_1|\cdots|g_*a_{i-1}]\otimes
g.
\end{array}
$$
where $\epsilon=|[a_1|\cdots|a_{i-1}]||[a_{i}|\cdots|a_n]|$. The
following lemma now holds by a direct calculation (where the reduced
chain complex is used):

\begin{lemma}Let $(C\otimes\Omega(C)\otimes k[G],b_g,\tilde{B})$ be as above. Then:
\begin{itemize}\item[$(a)$] $\tilde B^2=0$ and $b_g\tilde B+\tilde B b_g=id-g_*$.
\item[$(b)$]
$\tilde B$ commutes with the $k[G]$-action.
\end{itemize}\end{lemma}

With this lemma, we may consider the $G$-equivariant complex
$$\big(C\otimes\Omega(C)\otimes k[G]\big)/G=\big(C\otimes\Omega(C)\otimes
k[G]\big)\otimes_{k[G]}k,$$ where $\tilde b$ and $\tilde B$ descends
to $b$ and $B$, with $b^2=0$, $B^2=0$ and $bB+Bb=0$.

All above definitions can be generalized to the complete
case. Namely, for a complete DG coalgebra $C$ with a group
$G$-action, we may consider the complete tensor product of $C$ with
its complete cobar construction,
$$(C\hat\otimes\hat\Omega(C)\otimes k[G], \tilde b, \tilde B),$$
where $\tilde b$ and $\tilde B$ are the extensions of the usual
boundary operator $\tilde b$ and $\tilde B$ to the completion. Also
we may consider the $G$-equivariant complex
$$\big(C\hat\otimes\hat\Omega(C)\otimes k[G]/G, b, B\big).$$
And therefore, for the DG open Frobenius-like algebra $(A_{\rm
c}(\tilde M),C(\tilde M),\iota)$ on $\tilde M$, $A_{\rm c}(\tilde
M)\hat\otimes\hat\Omega(C(\tilde M))\otimes g$ and $C(\tilde
M)\hat\otimes\hat\Omega(C(\tilde M))\otimes g$ models the chain
complex of $L_g\tilde M$, and by Lemma \ref{equivchar} the
$G$-equivariant complex $\big(A_{\rm c}(\tilde
M)\hat\otimes\hat\Omega(C(\tilde M))\otimes k[G]\big)/G$ and
$\big(C(\tilde M)\hat\otimes\hat\Omega(C(\tilde M))\otimes
k[G]\big)/G$ models the chain complex of the free loop space
$LM$.

To simplify the notations we write $A_{\rm c}(\tilde
M)\hat\otimes\hat\Omega(C(\tilde M))\otimes g$ as $C_*(L_g\tilde
M)$, and $\big(A_{\rm c}(\tilde M)\hat\otimes\hat\Omega(C(\tilde M))
\otimes\Q[G]\big)/G$ as $C_*^G(\coprod L_g\tilde M)$ for short.

The loop product $\bullet$ of Chas and Sullivan is modeled as
follows:

\begin{definition}[Loop product]\label{lpofgeneral}
Let $(A_{\rm c}(\tilde M), C(\tilde M),\iota)$ be the DG
Frobenius-like algebra of $\tilde M$. Define a binary operator
$\tilde\bullet$ on $C_*(\coprod L_g\tilde M)$ as follows: for any
$$\a=x\otimes[a_1|\cdots|a_n]\otimes g\in C_*(L_g\tilde M)$$ and
$$\b=y\otimes[b_1|\cdots|b_m]\otimes h\in C_*(L_h\tilde M),$$ let $$
\a\tilde\bullet\b:=(-1)^{(|y|+1)|[a_1|\cdots|a_n]|}x\cdot
g_*^{-1}y\otimes[a_1|\cdots|a_n|b_1|\cdots|b_m]\otimes gh.
$$
On the $G$-equivariant chain complex $C_*^G(\coprod L_g\tilde M)$,
define a binary operator $\bullet$ as follows: for $[\a],[\b]\in
C_*^G(\coprod L_g\tilde M)$,
$$[\a]\bullet[\b]:=\Big[\a\tilde\bullet\sum_{g\in G}g_*\b\Big].$$
\end{definition}

\begin{lemma}The operator $\bullet$ does not depend on the choice of
the representatives and is well defined. Moreover, it commutes with
the boundary operator $b$.
\end{lemma}

\begin{proof}The fact that $\bullet$ commutes with $b$ follows from a direct computation
(compare with Definition \ref{modelofcs} in the simply connected case).
To show $\bullet$ does not depend on the choice of representatives,
take arbitrary $h,k\in G$,
$$[h_*\a]\bullet[k_*\b]=\Big[h_*\a\tilde\bullet\sum_{g\in G} g_*k_*\b\Big]
=\Big[h_*\a\tilde\bullet\sum_{g\in G}
g_*\b\Big]=\Big[h_*\a\tilde\bullet\sum_{g\in G}
h_*g_*\b\Big]=[h_*(\a\bullet\sum_{g\in G}g_* \b)]=[\a]\bullet[\b].
$$
Also since $\Q[G]$ acts on $C_*(\coprod L_g\tilde M)$ freely and
properly, and the differential forms are compactly supported,
$\bullet$ is well defined.
\end{proof}

Therefore we obtain a graded algebra on the homology of
$H_*(C_*^{G}(\coprod L_g\tilde M),b)$. As in the simply connected
case, such an algebra exactly models the loop product.

\begin{definition}[$*$ operator and the loop
bracket]\label{lbofgeneral} Let $(A_{\rm c}(\tilde M), C(\tilde
M),\iota)$ be the DG open Frobenius-like algebra of $\tilde M$.
Define a binary operator $\tilde*$ on $C_*(\coprod L_g\tilde M)$ as
follows: for any
$$\a=x\otimes[a_1|\cdots|a_n]\otimes g\in C_*(L_g\tilde M)$$ and
$$\b=y\otimes[b_1|\cdots|b_m]\otimes h\in C_*(L_h\tilde M),$$ let
$$
\a\tilde*\b:=\sum_i(-1)^{|y|+(|y|-1)|[a_{i+1}|\cdots|a_n]|}\varepsilon(a_iy)x\otimes
[a_1|\cdots|a_{i-1}|b_1|\cdots|b_m|h_*a_{i+1}|\cdots|h_*a_n]\otimes
gh.
$$
On the $G$-equivariant chain complex $C_*^G(\coprod L_g\tilde M)$,
define a binary operator $\ast$ as follows: for $[\a],[\b]\in
C_*^G(\coprod L_g\tilde M)$,
$$[\a]\ast[\b]:=\Big[\a\tilde*\sum_{g\in G}g_*\b\Big].$$
\end{definition}

\begin{lemma}[Gerstenhaber algebra of the free loop space]
Let $M$ and $\tilde M$ be as above. \begin{enumerate}\item[$(1)$] On
$C_*(\coprod L_g\tilde M)$,
$$b(\a\tilde*\b)=b\a\tilde*\b+(-1)^{|\a|+1}\a\tilde* b\b
+(-1)^{|\a|}(\a\tilde\bullet\b-(-1)^{|\a||\b|}h_*(h_*^{-1}\b\tilde\bullet
\a)).$$

\item[$(2)$]On $C_*^G(\coprod L_g\tilde M)$, the operator $*$ does
not depend on the choice of the representatives and hence is well defined.
Moreover,
$$b(\a\ast\b)=b\a\ast\b+(-1)^{|\a|+1}\a\ast b\b+(-1)^{|\a|}(\a\bullet\b-(-1)^{|\a||\b|}\b\bullet\a),$$
which means $\bullet$ is graded commutative on the homology
$\HH_*^G\Big(\coprod L_g\tilde M\Big)$.

\item[$(3)$]
The commutator of $*$ forms a degree one Lie algebra, which is
compatible with $\bullet$, making
$$\Big(\mathbb H_*^G\Big(\coprod L_g\tilde
M\Big),\bullet, \{,\}\Big)$$ a Gerstenhaber algebra.
\end{enumerate}
\end{lemma}

\begin{proof}
These results follow from a direct computation.
\end{proof}

\begin{theorem}[Batalin-Vilkovisky algebra]
Let $M$ be a smooth manifold and $\tilde M$ be its universal
covering. The homology
$$\Big(\mathbb H_*^G\Big(\coprod L_g\tilde M\Big),\bullet,B\Big)$$
forms a Batalin-Vilkovisky algebra, which coincides with the one
given by \cite{CS99}.
\end{theorem}

\begin{proof}
As in the above Definitions \ref{lpofgeneral} and \ref{lbofgeneral},
the homotopy operator defined in Lemma \ref{bvgersten} can be
applied here, which implies the theorem.
\end{proof}

The construction of the gravity algebra on the equivariant homology is
similar, and is left to the interested reader.

\renewcommand{\theequation}{a\arabic{equation}}
  % redefine the command that creates the equation no.
\setcounter{equation}{0}  % reset counter
\section*{Appendix: Proof of Lemmas \ref{lemma4.5} and \ref{bvgersten}}

\begin{proof}[Proof of Lemma \ref{lemma4.5}]
(1) comes immediately from the definitions of
$\bullet$ and $*$. We prove (2). In fact, up to sign,
\begin{eqnarray}&&
\gamma*(\a\bullet\b)-(\gamma*\a)\bullet\b-\a\bullet(\gamma*\b)\nonumber\\
&=&\sum_{i}\varepsilon(c_ixy)z\otimes
[c_1|\cdots|c_{i-1}|a_1|\cdots|a_n|b_1|\cdots|b_m|c_{i+1}|\cdots|c_l]\label{blue1}\\
&+&\sum_i\varepsilon(c_ix)yz\otimes
[c_1|\cdots|c_{i-1}|a_1|\cdots|a_n|c_{i+1}|\cdots|c_l|b_1|\cdots|b_m]\label{blue2}\\
&+&\sum_i
\varepsilon(c_iy)xz\otimes[a_1|\cdots|a_n|c_1|\cdots|c_{i-1}
|b_1|\cdots|b_m|c_{i+1}|\cdots|c_l],\label{blue3}
\end{eqnarray}
while
\begin{eqnarray}
&&b\circ h(\a,\b,\gamma)\nonumber\\
&=&\sum_{i<j}\varepsilon(c_ix)\varepsilon(c_jy)dz\otimes
[c_1|\cdots|c_{i-1}|a_1|\cdots|a_n|c_{i+1}|\cdots|c_{j-1}
|b_1|\cdots|b_m|c_{j+1}|\cdots|c_l]\label{g1}\\
&+&\sum_{i<j,r}\varepsilon(c_ix)\varepsilon(c_jy)z\otimes
[c_1|\cdots|c_{i-1}|a_1|\cdots|a_n|\cdots|dc_r|\cdots|
c_{j-1}|b_1|\cdots|b_m|\cdots|c_l]\label{g2}\\
&+&\sum_{i<j,r}\varepsilon(c_ix)\varepsilon(c_jy)z\otimes
[c_1|\cdots|c_{i-1}|a_1|\cdots|a_n|\cdots|c_r'|c_r''|\cdots|c_{j-1}
|b_1|\cdots|b_m|\cdots|c_l]\label{g3}\\
&+&\sum_{i<j,p}\varepsilon(c_ix)\varepsilon(c_jy)z\otimes
[c_1|\cdots|c_{i-1}|a_1|\cdots|da_p|\cdots|a_n|\cdots|c_{j-1}|
b_1|\cdots|b_m|\cdots|c_l]\label{g4}\\
&+&\sum_{i<j,p}\varepsilon(c_ix)\varepsilon(c_jy)z\otimes
[c_1|\cdots|c_{i-1}|a_1|\cdots|a_p'|a_p''|\cdots|a_n|\cdots|c_{j-1}
|b_1|\cdots|b_m|\cdots|c_l]\label{g5}\\
&+&\sum_{i<j,q}\varepsilon(c_ix)\varepsilon(c_jy)z\otimes
[c_1|\cdots|c_{i-1}|a_1|\cdots|a_n|\cdots|c_{j-1}|b_1|
\cdots|db_q|\cdots |b_m|\cdots|c_l]\label{g6}\\
&+&\sum_{i<j,q}\varepsilon(c_ix)\varepsilon(c_jy)z\otimes
[c_1|\cdots|c_{i-1}|a_1|\cdots|a_n|\cdots|c_{j-1}|b_1|
\cdots|b_q'|b_q''|\cdots |b_m|\cdots|c_l]\label{g7}\\
&+&\sum_{i<j}\varepsilon(c_ix)\varepsilon(c_jy)z'\otimes
[z''|c_1|\cdots|c_{i-1}|a_1|\cdots|a_n|\cdots|c_{j-1}
|b_1|\cdots|b_m|\cdots|c_l]\label{g8}\\
&+&\sum_{i<j}\varepsilon(c_ix)\varepsilon(c_jy)z'\otimes
[c_1|\cdots|c_{i-1}|a_1|\cdots|a_n|\cdots|c_{j-1}
|b_1|\cdots|b_m|\cdots|c_l|z''],\label{g9}
\end{eqnarray}and
\begin{eqnarray}
&&h(b\a,\b,\gamma)\nonumber\\
&=&\sum_{i<j}\varepsilon(c_idx)\varepsilon(c_jy)z\otimes
[c_1|\cdots|c_{i-1}|a_1|\cdots|a_n|c_{i+1}|\cdots|c_{j-1}
|b_1|\cdots|b_m|c_{j+1}|\cdots|c_l]\label{o1}\\
&+&\sum_{i<j,p}\varepsilon(c_ix)\varepsilon(c_jy)z\otimes
[c_1|\cdots|c_{i-1}|a_1|\cdots|da_p|\cdots|a_n|\cdots|c_{j-1}
|b_1|\cdots|b_m|\cdots|c_l]\label{o2}\\
&+&\sum_{i<j,p}\varepsilon(c_ix)\varepsilon(c_jy)z\otimes
[c_1|\cdots|c_{i-1}|a_1|\cdots|a_p'|a_p''|\cdots|a_n|\cdots|c_{j-1}
|b_1|\cdots|b_m|\cdots|c_l]\label{o3}\\
&+&\sum_{i<j}\varepsilon(c_ix')\varepsilon(c_jy)z
\otimes[c_1|\cdots|c_{i-1}|x''|a_1|\cdots|a_n|\cdots|c_{j-1}|b_1|\cdots|b_m|\cdots|c_l]
\label{o4}\\
&+&\sum_{i<j}\varepsilon(c_ix')\varepsilon(c_jy)z
\otimes[c_1|\cdots|c_{i-1}|a_1|\cdots|a_n|x''|\cdots|c_{j-1}|b_1|\cdots|b_m|\cdots|c_l],
\label{o5}
\end{eqnarray}and
\begin{eqnarray}
&&h(\a,b\b,\gamma)\nonumber\\
&=&\sum_{i<j}\varepsilon(c_ix)\varepsilon(c_jdy)z\otimes[c_1|\cdots|
c_{i-1}|a_1|\cdots|a_n|c_{i+1}|\cdots|c_{j-1}|b_1|\cdots|b_m|
c_{j+1}|\cdots|c_l]\label{o6}\\
&+&\sum_{i<j,q}\varepsilon(c_ix)\varepsilon(c_jy)z\otimes[c_1|\cdots|
c_{i-1}|a_1|\cdots|a_n|\cdots|c_{j-1}|b_1|\cdots|db_q|\cdots|b_m
|\cdots|c_l]\label{o7}\\
&+&\sum_{i<j,q}\varepsilon(c_ix)\varepsilon(c_jy)z\otimes[c_1|\cdots|
c_{i-1}|a_1|\cdots|a_n|\cdots|c_{j-1}|b_1|\cdots|b_q'|b_q''|\cdots|b_m
|\cdots|c_l]\label{o8}\\
&+&\sum_{i<j}\varepsilon(c_ix)\varepsilon(c_jy')z\otimes[c_1|\cdots|
c_{i-1}|a_1|\cdots|a_n|\cdots|c_{j-1}|y''|b_1|\cdots|b_m
|\cdots|c_l]\label{o9}\\
&+&\sum_{i<j}\varepsilon(c_ix)\varepsilon(c_jy')z\otimes[c_1|\cdots|
c_{i-1}|a_1|\cdots|a_n|\cdots|c_{j-1}|b_1|\cdots|b_m|y''
|\cdots|c_l],\label{o10}
\end{eqnarray}
and one has similar terms for $h(\a,\b,b\gamma)$. A straightforward
check shows that all terms cancel except for
(\ref{blue1})+(\ref{blue2})+(\ref{blue3}). Thus (2) is proved.
\end{proof}

\begin{proof}[Proof of Lemma \ref{bvgersten}]
First note that $A\hat\otimes\hat\Omega(C)$ embeds in
$C\hat\otimes\hat\Omega(C)$, so the operator $B$ is well defined.
For $\a=x\otimes[a_1|\cdots|a_n],\b=y\otimes[b_1|\cdots|b_m]\in
A\hat\otimes\hat \Omega(C)$, define
$$
\phi(\a,\b):=\sum_{i<j}\varepsilon(x)\varepsilon(a_jy)
a_i\otimes[a_{i+1}|\cdots|a_{j-1}|b_1|\cdots|b_m|a_{j+1}|\cdots|a_n|a_1|\cdots|a_{i-1}]
$$
and
$$\psi(\a,\b):=\sum_{k<l}\varepsilon(y)\varepsilon(b_lx) b_k\otimes
[b_{k+1}|\cdots|b_{l-1}|a_1|\cdots|a_n|b_{l+1}|\cdots|b_m|b_1|\cdots|b_{k-1}],
$$
and let $h=\phi+\psi$. We will show that $h$ thus defined satisfies
(\ref{bvhtpy}). Since $A\hat\otimes\hat\Omega(C)$ and
$C\hat\otimes\hat\Omega(C)$ have the same homology, (\ref{bv})
follows from (\ref{bvhtpy}).

In fact,
$\{\a,\b\}-(-1)^{|\a|}B(\a\bullet\b)-(-1)^{(|\b|+1)(|\a|+1)}\b\bullet
B(\a) +\a\bullet B(\b)$ contains two parts:
\begin{eqnarray}&&\sum_i\varepsilon(xy)a_i\otimes
[a_{i+1}|\cdots|a_n|b_1|\cdots|b_m|a_1|\cdots|a_{i-1}]\label{red4}\\
&+&\sum_i\varepsilon(a_iy)x\otimes[a_1|\cdots|a_{i-1}|b_1|\cdots|b_m|a_{i+1}|\cdots|a_n]\label{red1}\\
&+&\sum_i\varepsilon(x)a_iy\otimes[b_1|\cdots|b_m|a_{i+1}|\cdots|a_n|a_1|\cdots|a_{i-1}]\label{red6}
\end{eqnarray}
and
\begin{eqnarray}
&&\sum_k\varepsilon(xy)b_k\otimes[b_{k+1}|\cdots|b_m|a_1|\cdots|a_n|b_1|\cdots|b_{k-1}]\label{red2}\\
&+&\sum_k\varepsilon(b_kx)y\otimes[b_1|\cdots|b_{k-1}|a_1|\cdots|a_n|b_{k+1}|\cdots|b_m]\label{red3}\\
&+&\sum_k\varepsilon(y)b_kx\otimes[a_1|\cdots|a_n|b_{k+1}|\cdots|b_m|b_1|\cdots|b_{k-1}].\label{red5}
\end{eqnarray}
while
\begin{eqnarray}
&&b\phi(\a,\b)\nonumber\\
&=&\sum_{i<j}\varepsilon(x)\varepsilon(a_jy)da_i\otimes
[a_{i+1}|\cdots|a_{j-1}|b_1|\cdots|b_m|a_{j+1}|\cdots
|a_n|a_1|\cdots|a_{i-1}]\label{b2}\\
&+&\sum_{i<j,p}\varepsilon(x)\varepsilon(a_jy)
a_i\otimes[a_{i+1}|\cdots
|da_p|\cdots|a_{j-1}|b_1|\cdots|b_m|a_{j+1}|\cdots|a_n|
a_1|\cdots|a_{i-1}]\ \label{b1}\\
&+&\sum_{i<j,p}\varepsilon(x)\varepsilon(a_jy)
a_i\otimes[a_{i+1}|\cdots
|a_p'|a_p''|\cdots|a_{j-1}|b_1|\cdots|b_m|a_{j+1}|\cdots|a_n|
a_1|\cdots|a_{i-1}]\quad \label{b4}\\
&+&\sum_{i<j,q}\varepsilon(x)\varepsilon(a_jy)a_i\otimes
[a_{i+1}|\cdots|a_{j-1}|b_1|\cdots|db_q|\cdots|b_m|a_{j+1}|\cdots
|a_n|a_1|\cdots|a_{i-1}]\ \label{r9}\\
&+&\sum_{i<j,q}\varepsilon(x)\varepsilon(a_jy)a_i\otimes
[a_{i+1}|\cdots|a_{j-1}|b_1|\cdots|b_q'|b_q''|\cdots|b_m|a_{j+1}|\cdots
|a_n|a_1|\cdots|a_{i-1}]\label{r10}\\
&+&\sum_{i<j}\varepsilon(x)\varepsilon(a_jy)a_i'\otimes
[a_i''|a_{i+1}|\cdots|a_{j-1}|b_1|\cdots|b_m|
a_{j+1}|\cdots|a_n|a_1|\cdots|a_{i-1}]\label{b5}\\
&+&\sum_{i<j}\varepsilon(x)\varepsilon(a_jy)a_i''\otimes
[a_{i+1}|\cdots|a_{j-1}|b_1|\cdots|b_m|
a_{j+1}|\cdots|a_n|a_1|\cdots|a_{i-1}|a_i'].\label{b6}
\end{eqnarray}
and
\begin{eqnarray}&&\phi(b\a,\b)\nonumber\\
&=&\sum_{i<j,p}\varepsilon(a_jy)\varepsilon(x)a_i\otimes
[a_{i+1}|\cdots|da_p|\cdots|a_{j-1}|b_1|\cdots|b_m|a_{j+1}|\cdots|a_n|a_1|\cdots|a_{i-1}]\label{r1}\\
&+&\sum_{i<j,p}\varepsilon(a_jy)\varepsilon(x)a_i\otimes
[a_{i+1}|\cdots|a_p'|a_p''|\cdots|a_{j-1}|b_1|\cdots|b_m|a_{j+1}|\cdots|a_n|a_1|\cdots|a_{i-1}]\quad\label{r4}\\
&+&\sum_{i<j}\varepsilon(x)\varepsilon(a_jy)da_i\otimes
[a_{i+1}|\cdots|a_{j-1}|b_1|\cdots|b_m|a_{j+1}|\cdots|a_n|a_1|\cdots|a_{i-1}]\label{r2}\\
&+&\sum_{i<j}\varepsilon(x)\varepsilon(da_jy)a_i\otimes
[a_{i+1}|\cdots|a_{j-1}|b_1|\cdots|b_m|a_{j+1}|\cdots|a_n|a_1|\cdots|a_{i-1}]\label{r3}\\
&+&\sum_{i<j}\varepsilon(x)\varepsilon(a_j'y)a_i\otimes
[a_{i+1}|\cdots|a_{j-1}|b_1|\cdots|b_m|a_j''|a_{j+1}|\cdots|a_n|a_1|\cdots|a_{i-1}]\label{r7}\\
&+&\sum_{i<j}\varepsilon(x)\varepsilon(a_j''y)a_i\otimes
[a_{i+1}|\cdots|a_{j-1}|a_j'|b_1|\cdots|b_m|a_{j+1}|\cdots|a_n|a_1|\cdots|a_{i-1}]\label{r8}\\
&+&\sum_{i<j}\varepsilon(x)\varepsilon(a_jy)a_i'\otimes
[a_i''|\cdots|a_{j-1}|b_1|\cdots|b_m|a_{j+1}|\cdots|a_n|a_1|\cdots|a_{i-1}]\label{r5}\\
&+&\sum_{i<j}\varepsilon(x)\varepsilon(a_jy)a_i''\otimes
[a_{i+1}|\cdots|a_{j-1}|b_1|\cdots|b_m|a_{j+1}|\cdots|a_n|a_1|\cdots|a_{i-1}|a_i']\label{r6}\\
&+&\sum_i\varepsilon(x)\varepsilon(a_i''y)a_i'\otimes
[b_1|\cdots|b_m|a_{i+1}|\cdots|a_n|a_1|\cdots|a_{i-1}]\label{black6}\\
&+&\sum_i\varepsilon(x')\varepsilon(a_iy)x''\otimes
[a_1|\cdots|a_{i-1}|b_1|\cdots|b_m|a_{i+1}|\cdots|a_n]\label{black1}\\
&+&\sum_i\varepsilon(x')\varepsilon(x''y)a_i\otimes[a_{i+1}|\cdots|a_n|b_1|\cdots|b_m|a_1|\cdots|a_{i-1}],\label{black2}
\end{eqnarray}
and
\begin{eqnarray}
&&\phi(\a,b\b)\nonumber\\
&=&\sum_{i<j}\varepsilon(x)\varepsilon(a_jdy)a_i\otimes
[a_{i+1}|\cdots|a_{j-1}|b_1|\cdots|b_m|a_{j+1}|\cdots|a_n|a_1|\cdots|a_{i-1}]\label{b3}\\
&+&\sum_{i<j,q}\varepsilon(x)\varepsilon(a_jy)a_i\otimes
[a_{i+1}|\cdots|a_{j-1}|b_1|\cdots|db_q|\cdots|b_m|
a_{j+1}|\cdots|a_n|a_1|\cdots|a_{i-1}]\label{b9}\\
&+&\sum_{i<j,q}\varepsilon(x)\varepsilon(a_jy)a_i\otimes
[a_{i+1}|\cdots|a_{j-1}|b_1|\cdots|b_q'|b_q''|\cdots|b_m|
a_{j+1}|\cdots|a_n|a_1|\cdots|a_{i-1}]\label{b10}\\
&+&\sum_{i<j}\varepsilon(x)\varepsilon(a_jy')a_i\otimes
[a_{i+1}|\cdots|a_{j-1}|y''|b_1|\cdots|b_m|a_{j+1}|\cdots
|a_n|a_1|\cdots|a_{i-1}]\label{b8}\\
&+&\sum_{i<j}\varepsilon(x)\varepsilon(a_jy')a_i\otimes
[a_{i+1}|\cdots|a_{j-1}|b_1|\cdots|b_m|y''|a_{j+1}|\cdots
|a_n|a_1|\cdots|a_{i-1}].\label{b7}
\end{eqnarray}
Note that (\ref{r1}) and (\ref{b1}) are identical, so are (\ref{r2})
and (\ref{b2}), (\ref{r3}) and (\ref{b3}), (\ref{r4}) and
(\ref{b4}), (\ref{r5}) and (\ref{b5}), (\ref{r6}) and (\ref{b6}),
(\ref{r7}) and (\ref{b7}), (\ref{r8}) and (\ref{b8}), (\ref{r9}) and
(\ref{b9}), and (\ref{r10}) and (\ref{b10}), therefore the remaining
terms of $b\phi(\a,\b)-\phi(b\a,\b)-\phi(\a,b\b)$ are exactly
(\ref{red4})+(\ref{red1})+(\ref{red6}).

Similarly $b\psi(\a,\b)-\psi(b\a,\b)-\psi(\a,b\b)$ is equal to
(\ref{red2})+(\ref{red3})+(\ref{red5}). \end{proof}

\end{document}